\documentclass{amsart}

\input{diagrams}
\usepackage{amsmath,amsthm,verbatim,amssymb,amsfonts,amscd,graphics}
\newtheorem{thm}{Theorem}[section] 

\newtheorem{cor}[thm]{Corollary}
\newtheorem{defn}[thm]{Definition}
\newtheorem{example}[thm]{Example}

\newtheorem{lemma}[thm]{Lemma}
\newtheorem{prop}[thm]{Proposition}
\newtheorem{remark}[thm]{Remark}
\def\abs#1{|#1|}
\newarrow{ul}---->
\newarrow{Backwards}<----

\begin{document}
\title[Alexander Modules of Hypersurface Complements]{Intersection Homology and Alexander Modules of Hypersurface Complements}
\author{Laurentiu Maxim}
\date{December, 29}
\maketitle

\begin{abstract}
Let $V$ be a degree $d$, reduced hypersurface in $\mathbb{CP}^{n+1}$, $n \geq 1$,
and fix a generic hyperplane, $H$. Denote by $\mathcal{U}$  the
(affine) hypersurface complement, $\mathbb{CP}^{n+1}- V \cup H$,
and let $\mathcal{U}^c$ be the infinite cyclic covering of
$\mathcal{U}$ corresponding to the kernel of the linking number
homomorphism. Using intersection homology theory, we give a new
construction of the  Alexander modules
$H_i(\mathcal{U}^c;\mathbb{Q})$ of the hypersurface complement and
show that, if $i \leq n$, these  are
torsion over the ring of rational Laurent polynomials.
We also obtain obstructions on the associated global polynomials.
Their zeros are roots of unity of order $d$ and are entirely determined
by the local topological information encoded by the link pairs of singular strata
of a stratification of the pair
$(\mathbb{CP}^{n+1},V)$. As an application, we give obstructions
on the eigenvalues of monodromy operators associated to
 the Milnor fibre of a projective hypersurface arrangement.

\end{abstract}

\tableofcontents

\section{Introduction}
Intersection homology is the correct theory to extend results from manifolds to singular varieties:
e.g., Morse theory, Lefschetz (weak and hard) theorems, Hodge decompositions,
but especially Poincar\'{e} duality (which motivates the theory).
It is therefore natural to use it in order to describe
topological invariants associated with algebraic varieties.

We will use intersection homology for the study of Alexander modules of hypersurface complements.
These are invariants introduced and studied by Libgober in a sequence of papers
\cite{Li1}, \cite{Li}, \cite{Li3}, \cite{Li4}, but see also \cite{Di3}. Their definition goes as follows.
Let $V$ be a degree $d$, reduced, projective hypersurface in $\mathbb{CP}^{n+1}$, $n \geq 1$; let $H$ be a
fixed hyperplane which we call 'the hyperplane at infinity'; set $\mathcal{U}:=\mathbb{CP}^{n+1}- V \cup H$.
(Alternatively, let $X \subset \mathbb{C}^{n+1}$ be a reduced affine hypersurface and
$\mathcal{U}:=\mathbb{C}^{n+1}- X$.)
Then $H_1(\mathcal{U}) \cong \mathbb{Z}^s$, where $s$ is the number of components of $V$,
and one proceeds, as in  classical knot theory, to define Alexander-type invariants of the hypersurface $V$. More
precisely, the rational homology groups of any infinite cyclic cover of $\mathcal{U}$ become, under the action of the
group of covering transformations,
modules over the ring $\Gamma$ of rational Laurent polynomials. The modules $H_i(\mathcal{U}^c;\mathbb{Q})$ associated to the infinite cyclic cover $\mathcal{U}^c$ defined by the kernel of the total linking number homomorphism, are called \emph{the Alexander modules of the hypersurface
complement}. Note that, since $\mathcal{U}$ has the homotopy type of a finite CW complex of
dimension $\leq n+1$, the Alexander modules $H_{i}({\mathcal{U}}^{c};\mathbb{Q})$ are trivial for $i > n+1$ and
$H_{n+1}({\mathcal{U}}^{c};\mathbb{Q})$ is $\Gamma$-free. Thus, of a particular interest are
the Alexander modules $H_i(\mathcal{U}^c;\mathbb{Q})$ for $i \leq n$.

Libgober showed (\cite{Li}, \cite{Li3}, \cite{Li4}, \cite{Li5})
that if $V$ has only isolated singularities (including at
infinity), then $\tilde{H}_{i}({\mathcal{U}}^{c};\mathbb{Z})=0$
for $i <n$, and $H_{n}({\mathcal{U}}^{c};\mathbb{Q})$ is a torsion
$\Gamma$-module. Also, if $\delta_n (t)$ denotes the polynomial
associated to the torsion module
$H_{n}({\mathcal{U}}^{c};\mathbb{Q})$, then $\delta_n (t)$ divides
(up to a power of $t-1$) the product of the local Alexander
polynomials of the algebraic knots around the isolated singular
points. If  $H$ is generic (hence $V$ has no singularities at
infinity), then the zeros of $\delta_n (t)$ are roots of unity of
order $d$, and  $H_{n}({\mathcal{U}}^{c};\mathbb{Q})$ is a
semi-simple module annihilated by $t^d-1$.\newline

The aim of this paper is to provide generalizations of these
results to the case of hypersurfaces with non-isolated
singularities. We will assume that $H$ is generic, i.e.,
transversal to all strata of a Whitney stratification of $V$. Using
intersection homology theory, we will give a new description of
the Alexander modules of the hypersurface complement. These will
be realized as
 intersection homology groups of $\mathbb{CP}^{n+1}$, with a certain local coefficient
system with stalk $\Gamma:=\mathbb{Q}[t,t^{-1}]$, defined on
$\mathcal{U}$.  Therefore, we will
have at our disposal the apparatus of intersection homology and
derived categories to study the Alexander modules of the
complement.\newline

The structure of the paper is as follows:

In Section 2, we recall the definitions and main properties of the
Alexander modules of the hypersurface complement,
$H_i(\mathcal{U}^c;\mathbb{Q})$. We also show that, if $V$ is in
general position at infinity, has no codimension one
singularities, and is a rational homology
manifold, then for $i \leq n$, the modules
$H_i(\mathcal{U}^c;\mathbb{Q})$ are torsion and their associated
polynomials do not contain factors $t-1$ (see Proposition 2.1).

In Section 3, we realize the Alexander modules of complements to hypersurfaces in general position at infinity as
intersection homology modules.
In Section 3.1, we first construct the intersection Alexander modules of the hypersurface $V$. More precisely,
by choosing a Whitney stratification $\mathcal{S}$ of $V$ and a generic hyperplane, $H$,
we obtain a stratification of the pair $(\mathbb{CP}^{n+1},V \cup H)$. We  define a
local system  $\mathcal{L}_{H}$ on $\mathcal{U}$,
with stalk $\Gamma := \mathbb{Q} [t, t^{-1}]$ and action by an element
$\alpha\in\pi_1(\mathcal{U})$ determined by multiplication by
$t^{\text{lk}(V \cup -dH,\alpha)}$.
Then, for any perversity ${\bar p}$, the intersection homology complex
$IC_{\bar p}^{\bullet}:=IC_{\bar p} ^{\bullet} (\mathbb{CP}^{n+1}, \mathcal{L}_{H})$ is defined  by using  Deligne's axiomatic construction (\cite{B}, \cite{GM2}). The
\emph{intersection Alexander modules of the hypersurface $V$} are then defined as  hypercohomology groups of the
middle-perversity intersection homology complex:
$IH_i^{\bar m}(\mathbb{CP}^{n+1};\mathcal{L}_{H}):=
\mathbb{H}^{-i}(\mathbb{CP}^{n+1};IC_{\bar m}^{\bullet})$.

In Section 3.2, we first prove the key technical lemma, which
asserts that the restriction to $V \cup H$ of the intersection
homology complex $IC_{\bar m} ^{\bullet}$ is quasi-isomorphic to the zero complex (see
Lemma 3.1). As a corollary, it follows that \emph{the intersection
Alexander modules of $V$ coincide with the Alexander modules of
the hypersurface complement}, i.e. there is an isomorphism of
$\Gamma$-modules: $IH_{*} ^{\bar m}
(\mathbb{CP}^{n+1};\mathcal{L}_{H}) \cong H_{*}(\mathcal{U}
^{c};\mathbb{Q})$. From now on, we will study the intersection
Alexander modules in order to obtain results on the Alexander
modules of the complement. Using the superduality isomorphism for
the local finite type codimension two embedding $V \cup H \subset
\mathbb{CP}^{n+1}$, and the peripheral complex associated with the
embedding (see \cite{ShCa}), we show that the $\Gamma$-modules
$IH_{i} ^{\bar m} (\mathbb{CP}^{n+1};\mathcal{L}_{H})$ are
\emph{torsion} if $i \leq n$ (see Corollary 3.8). Therefore the
classical Alexander modules of the hypersurface complement are
torsion in the range $i \leq n$. We denote their associated
polynomials by $\delta_i(t)$ and call them \emph{the global
Alexander polynomials of the hypersurface}.

Section 4 contains the main theorems of the paper, which provide obstructions on the prime divisors of the  polynomials $\delta_i(t)$.
Our results are extensions to the case of
hypersurfaces with general singularities of the results proven by A. Libgober for hypersurfaces with isolated singularities
(\cite{Li}, \cite{Li3}, \cite{Li4}).

  The first theorem
gives a characterization of the zeros of the global polynomials and generalizes Corollary $4.8$ of \cite{Li}:
\begin{thm}(see Theorem 4.1)\newline
If $V$ is an n-dimensional reduced projective hypersurface of degree $d$, transversal to the hyperplane at infinity,
then  the zeros of the global Alexander polynomials $\delta_{i}(t)$, $i \leq n$, are roots of unity of order $d$.
\end{thm}

The underlying idea of this paper is to use local topological information associated with a singularity to
describe some global topological invariants of algebraic varieties. We provide a general divisibility
result which
restricts the prime factors of the global Alexander polynomial $\delta_i(t)$  to those
of the local  Alexander polynomials of the link knots around the singular strata. More precisely, we prove the following:
\begin{thm}(see Theorem 4.2)\newline
Let $V$ be a reduced hypersurface in $\mathbb{CP}^{n+1}$, which is
transversal to the hyperplane at infinity, $H$. Fix an arbitrary
irreducible component of $V$, say $V_1$. Let $\mathcal{S}$ be a
stratification of the pair $(\mathbb{CP}^{n+1},V)$. Then for a
fixed integer $i$ ($1 \leq i \leq n$), the prime factors of the
global Alexander polynomial $\delta_i(t)$ of $V$ are among the
prime factors of local polynomials $\xi^s _{l}(t)$ associated to
the local Alexander modules $H_l (S^{2n-2s+1}-K^{2n-2s-1};\Gamma)$
of
 link pairs $(S^{2n-2s+1},K^{2n-2s-1})$
of components of strata $S \in \mathcal{S}$ such that: $S \subset
V_1$, $n-i \leq s=\text{dim}{S} \leq n$, and $l$ is in the range
$2n-2s-i \leq l \leq n-s$.
\end{thm}
For hypersurfaces with isolated singularities, the above theorem
can be strengthened to obtain a result similar to Theorem 4.3 of
\cite{Li} (in the absence of singular points at infinity):
\begin{thm}(see Theorem 4.5)\newline
Let $V$ be a reduced hypersurface in $\mathbb{CP}^{n+1}$ ($n
\geq 1$), which is transversal to the hyperplane at infinity, H,
and has only isolated singularities. Fix an irreducible component
of $V$, say $V_1$. Then  $\delta_n(t)$ divides (up to a power of $(t-1)$) the product $\prod_{p \in
V_1 \cap \text{Sing}(V)} \Delta_{p}(t)$ of the local Alexander
polynomials of links of the singular points $p$ of $V$ which are
contained in $V_1$.
\end{thm}
We end the section  by relating the intersection Alexander modules
of $V$ to the modules 'at infinity'. We prove the following extension of Theorem $4.5$ of \cite{Li}:
\begin{thm}(see Theorem 4.7)\newline
Let $V$ be a reduced hypersurface of degree $d$ in $\mathbb{CP}^{n+1}$, which is transversal to the
hyperplane at infinity, $H$.
Let $S_{\infty}$ be a sphere of sufficiently large radius in $\mathbb{C}^{n+1}=
\mathbb{CP}^{n+1} -H$. Then for all $i < n$,
$$IH_{i} ^{\bar m} (\mathbb{CP}^{n+1};\mathcal{L}_{H}) \cong \mathbb{H}^{-i-1}(S_{\infty};IC^{\bullet} _{\bar m}) \cong
H_{i} (\mathcal{U}^c_{\infty};\mathbb{Q})$$ and $IH_{n} ^{\bar m}
(\mathbb{CP}^{n+1};\mathcal{L}_{H})$ is a quotient of
$\mathbb{H}^{-n-1}(S_{\infty};IC^{\bullet} _{\bar m}) \cong H_{n}
(\mathcal{U}^c_{\infty};\mathbb{Q})$, where
$\mathcal{U}^c_{\infty}$ is the infinite cyclic cover of
$S_{\infty}-(V \cap S_{\infty})$ corresponding to the linking
number with $V \cap S_{\infty}$.
\end{thm}
Based on a private communication of A. Libgober (\cite{Li6}), we
note that the above theorem has as a corollary the semi-simplicity
of the Alexander modules of the hypersurface complement (see
Proposition 4.9).

In Section 5, we apply the preceding results to the case of a
hypersurface $V \subset \mathbb{CP}^{n+1}$, which is a projective
cone over a reduced hypersurface $Y \subset \mathbb{CP}^n$. We
first note that Theorem 4.2 translates into divisibility results
for the characteristic polynomials of the monodromy operators
acting on the Milnor fiber $F$ of the projective arrangement
defined by $Y$ in $\mathbb{CP}^n$. We obtain the following result,
similar to those obtained by A. Dimca in the case of isolated
singularities (\cite{Di}, \cite{Di2}):
\begin{prop}(see Proposition 5.1)\newline
Let $Y=(Y_i)_{i=1,s}$ be a hypersurface arrangement in
$\mathbb{CP}^n$, and fix an arbitrary component, say $Y_1$.  Let
$F$ be the  Milnor fibre of the arrangement. Fix a Whitney
stratification of the pair $(\mathbb{CP}^n,Y)$ and denote by
$\mathcal{Y}$ the set of (open) singular strata. Then for $q \leq
n-1$, a prime $\gamma \in \Gamma$ divides the characteristic
polynomial $P_q(t)$ of the monodromy operator $h_q$ only if
$\gamma$ divides one of the polynomials $\xi^s _{l}(t)$ associated
to the local Alexander modules
$H_{l}(S^{2n-2s-1}-K^{2n-2s-3};\Gamma)$ corresponding to
 link pairs $(S^{2n-2s-1},K^{2n-2s-3})$
of components of strata $\mathcal{V} \in \mathcal{Y}$ of complex
dimension $s$ with $\mathcal{V} \subset Y_1$, such that: $n-q-1
\leq s \leq n-1$ and $2(n-1)-2s-q \leq l \leq n-s-1$.
\end{prop}
As a consequence, we obtain obstructions on the eigenvalues of the
monodromy operators (see Corollary 5.3), similar to those obtained
by Libgober in the case of hyperplane arrangements (\cite{Li100}),
or Dimca in the case of curve arrangements (\cite{Di2}).

Section 6 deals with examples. We show, by explicit calculations,
how to apply the above theorems in obtaining information on the
global Alexander polynomials of a hypersurface in general position
at infinity.\newline

\emph{Note.} Our overall approach makes use of sheaf theory and
the language of derived categories (\cite{B}, \cite{ShCa},
\cite{GM2}) and we rely heavily on the material in these
references. We will always use the indexing conventions of
\cite{GM2}.

{\bf Acknowledgements.} I would like to express my deep gratitude
to my advisor, Professor Julius Shaneson, for encouragement and
advice. I am grateful to Greg Friedman, Alexandru Dimca, Anatoly
Libgober, and Mark Andrea de Cataldo for useful discussions. This
work is part of the author's thesis.

\emph{Note.}  After reading this paper, A. Dimca has obtained
alternative proofs for some of the results, including Corollary
3.8, Theorem 4.1 and Theorem 4.2. Subsequently, A. Libgober has
found different proofs for Corrolary 3.8, Theorem 4.1 and Theorem
4.7.


\section{Alexander modules of hypersurface complements}
In this section we recall the definition and main results on the Alexander modules and polynomials of
hypersurface complements. We also consider the special case of hypersurfaces which are rational homology manifolds.


\subsection{Definitions} Let $X$ be a connected CW complex, and let $\pi_{X}: \pi_{1}(X) \to \mathbb{Z}$ be an epimorphism.
 We denote by
$X^{c}$ the $\mathbb{Z}$-cyclic covering associated to the kernel of the morphism $\pi_{X}$.
The group of covering transformations of $X^{c}$ is infinite cyclic and acts
 on $X^{c}$ by a covering homeomorphism $h$.
Thus, all the groups $H_{\ast}(X^{c}; A)$, $H^{\ast}(X^{c}; A)$ and $\pi_{j}(X^{c}) \otimes A$
for $j >1$ become in the usual way $\Gamma_{A}$-modules, where $\Gamma_{A}=A[t,t^{-1}]$, for any ring $A$.
These are called \emph{the Alexander modules of the pair $(X,\pi_{X})$}.

If $A$ is a field, then the ring $\Gamma_{A}$ is a PID. Hence
any torsion $\Gamma_{A}$-module $M$ of finite type has a well-defined associated order (see \cite{Mi2}).
This is called
\emph{the Alexander polynomial of the torsion $\Gamma_{A}$-module
$M$}  and denoted by $\delta_M(t)$. We regard the trivial module $(0)$ as a
torsion module whose associated polynomial is $\delta(t)=1$.

With these notations, we have the following simple fact: let $f :M \to N$ be an epimorphism of $R$-modules,
where $R$ is a PID and $M$ is torsion of finite type. Then $N$ is torsion
of finite type and $\delta_N(t)$ divides $\delta_M(t)$.


\subsection{Alexander modules of hypersurface complements}
To fix notations for the rest of the paper, let $V$ be a \emph{reduced} hypersurface  in $\mathbb{CP}^{n+1}$,
defined by a degree $d$ homogeneous equation: $f=f_1\cdots f_s =0$, where $f_i$ are the irreducible factors of $f$
and $V_i=\{f_i=0\}$ the irreducible components of $V$. We will assume that $V$ is \emph{in general position at infinity},
i.e. we choose a generic hyperplane $H$ (transversal to all singular strata in a stratification of $V$) which we call
'the hyperplane at infinity'.
Let $\mathcal{U}$ be the (affine) hypersurface complement: $\mathcal{U}=\mathbb{CP}^{n+1} - (V \cup H)$.
Then $H_1 (\mathcal{U}) \cong \mathbb{Z}^s$ (\cite{Di}, (4.1.3), (4.1.4)), generated by the meridian loops $\gamma_i$
about the non-singular part of each irreducible component $V_i$, $i=1,\cdots, s$.
If $\gamma_\infty$ denotes the meridian  about the hyperplane at infinity, then in $H_1(\mathcal{U})$ there is a relation:
$\gamma_\infty + \sum {d_i \gamma_i} = 0$, where $d_i=deg(V_i)$.
We consider the infinite cyclic cover $\mathcal{U}^c$ of $\mathcal{U}$ defined by the kernel of the total linking number
homomorphism $lk : \pi_1(\mathcal{U}) \to \mathbb{Z}$, which maps all the meridian generators to $1$,
and thus any loop  $\alpha$ to $\text{lk} (\alpha, V \cup -dH)$.
Note that $lk$ coincides with the homomorphism $\pi_1(\mathcal{U}) \to \pi_1(\mathbb{C}^*)$ induced by the
polynomial map defining the affine hypersurface $V_{aff}:= V - V\cap H$ (\cite{Di}, p. 76-77).
The \emph{Alexander modules of the hypersurface complement} are defined as
$H_{i}({\mathcal{U}}^{c};\mathbb{Q})$.

Since $\mathcal{U}$ has the homotopy type of a finite CW complex of dimension $\leq n+1$ (\cite{Di} (1.6.7), (1.6.8)),
it follows that all the associated Alexander modules are of finite type over $\Gamma_{A}$, but in general not over $A$.
It also follows that the Alexander modules $H_{i}({\mathcal{U}}^{c};\mathbb{Q})$ are trivial for $i > n+1$ and
$H_{n+1}({\mathcal{U}}^{c};\mathbb{Q})$ is free over $\Gamma_{\mathbb{Q}}$ (\cite{Di3}).

Note that if $V$ has no codimension one singularities (e.g. if $V$ is normal), then
the fundamental group of $\mathcal{U}$ is infinite cyclic (\cite{Li}, Lemma 1.5) and ${\mathcal{U}}^{c}$ is the universal
cover of $\mathcal{U}$.
In particular, this condition if satisfied if $n \geq 2$ and $V$ has only isolated singularities: if this is the case,
Libgober shows (\cite{Li}) that
$\tilde{H}_{i}({\mathcal{U}}^{c};\mathbb{Z})=0$ for $i <n$, and $H_{n}({\mathcal{U}}^{c};\mathbb{Q})$ is a torsion
$\Gamma_{\mathbb{Q}}$-module. Also, one has the isomorphisms of the Alexander $\Gamma_{\mathbb{Z}}$-modules:
$\pi_{n} (\mathcal{U}) = \pi_{n} ({\mathcal{U}}^{c}) = H_{n}({\mathcal{U}}^{c};\mathbb{Z})$;
if we denote by $\delta_n (t)$ the polynomial associated to the torsion module $H_{n}({\mathcal{U}}^{c};\mathbb{Q})$,
then Theorem
4.3 of \cite{Li} asserts that $\delta_n (t)$ divides the product $\prod_{i=1} ^s \Delta_{i}(t) \cdot (t-1)^k$
of the Alexander polynomials of links of the singular points of $V$. Moreover, the zeros of $\delta_n (t)$ are roots of
unity of order $d = deg(V)$ (\cite{Li}, Corollary 4.9).

\emph{Note.} Libgober's   divisibility theorem (\cite{Li}, Theorem 4.3) holds for hypersurfaces with isolated
singularities, \emph{including at infinity} (and $n \geq 1$). However, for non-generic $H$ and for  hypersurfaces with more general
singularities, the Alexander modules
$H_{i}({\mathcal{U}}^{c};\mathbb{Q})$ ($i \leq n$) are not torsion in general. Their $\Gamma_{\mathbb{Q}}$ rank is calculated in \cite{Di3},
Theorem 2.10(v). We will show that if $V$ is a reduced hypersurface, in general position
at infinity, then the modules $H_{i}({\mathcal{U}}^{c};\mathbb{Q})$  are torsion for $i \leq n$.


\subsection{Rational homology manifolds} Recall that a $n$-dimensional complex variety $V$ is called a
\emph{rational homology manifold}, or is said to be
\emph{rationally smooth}, if for all points $x \in V$ we have:

\begin{equation*}
H_{i}(V,V-x;\mathbb{Q}) \cong
\begin{cases}
\mathbb{Q} \ ,& i=2n\\
0, & i \neq 2n.
\end{cases}
\end{equation*}

A rational homology manifold of dimension $n$ has pure dimension $n$ as a complex variety. Rational homology manifolds
may be thought of as 'nonsingular for the purposes of rational homology'. For example, Poincar\'{e} and Lefschetz duality
hold for them in rational homology. The Lefschetz hyperplane section theorem also holds.
Examples of rational homology manifolds include
complex varieties having rational homology spheres as links of singular strata.

Note that, if $V$ is a projective hypersurface having rational homology spheres as links of singular strata
 and $H$ is a generic hyperplane, then $V \cap H$ is a rational homology manifold: indeed, by the transversality assumption,
 the link in $V \cap H$ of a
 stratum $S \cap H$ (for $S$ a stratum of $V$) is the same as the link in $V$
 of $S$.

 As a first example when the Alexander modules $H_{i}({\mathcal{U}}^{c};\mathbb{Q})$, $i \leq n$, are  torsion, we mention the
 following (compare \cite{Li}, Lemma 1.7, 1.12):

\begin{prop}
Let $V$ be a degree $d$ projective hypersurface in $\mathbb{CP}^{n+1}$, and let $H$ be a generic
hyperplane. Assume that $V$ has no
codimension one singularities and that $V$ is a  rational homology manifold.
Then for $i \leq n$, $H_{i}({\mathcal{U}}^{c};\mathbb{Q})$ is a  torsion $\Gamma_{\mathbb{Q}}$-module
and $\delta_i (1) \neq 0$, where $\delta_i (t)$ is the associated Alexander polynomial.
\end{prop}
\begin{proof}
Recall that, under our assumptions, ${\mathcal{U}}^{c}$ is the infinite cyclic and universal cover of $\mathcal{U} = \mathbb{CP}^{n+1} - (V \cup H)$.
We will use  Milnor's exact sequence (\cite{Li}, \cite{Di3}):
$$\cdots \to H_{i}({\mathcal{U}}^{c};\mathbb{Q}) \to H_{i}({\mathcal{U}}^{c};\mathbb{Q}) \to H_{i}(\mathcal{U};\mathbb{Q})
\to H_{i-1}({\mathcal{U}}^{c};\mathbb{Q}) \to \cdots $$ where the first morphism is  multiplication by $t-1$. We claim
that $H_{i}(\mathcal{U};\mathbb{Q}) \cong 0$ for $2 \leq i \leq n$, hence the multiplication by $t-1$ in
$H_{i}({\mathcal{U}}^{c};\mathbb{Q})$ is surjective ($2 \leq i \leq n$). Therefore its cyclic decomposition has neither free
summands nor summands of the form $\Gamma_{\mathbb{Q}}/(t-1)^{r} \Gamma_{\mathbb{Q}}$, with $r \in \mathbb{N}$. On the other
hand, $H_{1}({\mathcal{U}}^{c};\mathbb{Q}) \cong \pi_1(\mathcal{U}^{c})\otimes \mathbb{Q} \cong 0$.

 Suppose that $k$ is the
dimension of the singular locus of $V$. Then, by our assumptions, $n-k \geq 2$.
Let $L \cong \mathbb{CP}^{n-k}$ be a generic linear subspace. Then, by transversality, $L \cap V$ is a
non-singular hypersurface in $L$, transversal to the hyperplane at infinity, $L \cap H$. Therefore, by Corollary 1.2 of
\cite{Li}, $L \cap \mathcal{U}$ is homotopy equivalent to $S^{1} \vee S^{n-k} \vee \cdots \vee S^{n-k}$.
Thus, by  Lefschetz hyperplane section theorem (applied $k+1$ times) we obtain: $H_i
(\mathcal{U};\mathbb{Q}) \cong H_i (L \cap \mathcal{U};\mathbb{Q})=0 \ , \ 2 \leq i \leq n-k-1$.

For $n-k \leq i \leq n$ we have:
$H_i (\mathcal{U};\mathbb{Q}) \cong H_{i+1} (\mathbb{CP}^{n+1} - H,\mathbb{CP}^{n+1} - (V \cup H);\mathbb{Q})$, as follows from
the exact sequence of the pair $(\mathbb{CP}^{n+1} - H,\mathbb{CP}^{n+1} - (V \cup H))$. Using duality, one can
identify this with $H^{2n+1-i} (V \cup H, H;\mathbb{Q})$. And by excision, this group is isomorphic to
$H^{2n+1-i} (V , V \cap H;\mathbb{Q})$. Let $u$ and $v$ denote the inclusion of $V - V \cap H$ and respectively $V \cap H$
into $V$. Then the distinguished triangle $u_! u^! \mathbb{Q} \to \mathbb{Q} \to v_* v^* \mathbb{Q} \overset{[1]}{\to}$ (where we regard
$\mathbb{Q}$
as a constant sheaf on $V$),
upon applying the hypercohomology with compact support functor, yields the isomorphism:
$H^{2n+1-i} (V , V \cap H;\mathbb{Q}) \cong H_c^{2n+1-i} (V - V \cap H;\mathbb{Q})$ (see \cite{Di2}, Remark 2.4.5.(iii)).
By Poincar\'{e} duality
over $\mathbb{Q}$, the latter is isomorphic to $H_{i-1} (V - V \cap H;\mathbb{Q})$. The Lefschetz theorem on generic
hyperplane complements in hypersurfaces (\cite{Di4}, p. 476) implies that $V - V \cap H$ is homotopy equivalent to a wedge
of spheres $S^n$. Therefore, $H_{i-1} (V - V \cap H;\mathbb{Q}) \cong 0$ for $0 < i-1 <n$, i.e. for $2 \leq n-k \leq i \leq
n$. This finishes the proof of the proposition.
\end{proof}


\section{Intersection Homology and Alexander modules}

Using intersection homology theory, we will give a new
construction of the Alexander modules of complements of
hypersurfaces in general position at infinity. The advantage of
the new approach is the use of the powerful language of sheaf
theory and derived categories (\cite{GM2}, \cite{B}) in the study
of the Alexander invariants associated with singular
hypersurfaces. This will alow us to obtain generalizations to classical results known only
in the case of hypersurfaces with isolated singularities (\cite{Li}).
 For a quick introduction to derived
categories, the reader is advised to consult \cite{Ma}, $\S 1$.
When dealing with intersection homology, we will always use the
indexing conventions of \cite{GM2}.


\subsection{Intersection Alexander Modules}

(1) \ \ A knot is a sub-pseudomanifold $K^{n} \subset S^{n+2}$ of a sphere; it is said to be of \emph{finite (homological) type} if the
homology groups $H_i (S^{n+2}-K;\Gamma)$ with local coefficients in $\Gamma:=\mathbb{Q}[t,t^{-1}]$ are finite dimensional
over $\mathbb{Q}$. Here $\Gamma$ denotes the local system on $S^{n+2}-K$, with
stalk $\Gamma$, and it corresponds to the
representation $\alpha \mapsto t^{\text{lk}(K,\alpha)}$, $\alpha\in\pi_1(S^{n+2}-K)$, where $\text{lk}(K,\alpha)$ is the
linking number of $\alpha$ with $K$ (see \cite{ShCa}).

A sub-pseudomanifold $X$ of a manifold $Y$ is said to be of \emph{finite local type} if the link of each
component of any stratification of the pair $(Y,X)$ is of finite type. Note that a sub-pseudomanifold is of finite local type if and only if it
has one stratification with links of finite type. It is also not hard to see that the link pairs of components of strata of
a sub-pseudomanifold of finite local type also have finite local type (\cite{ShCa}). Algebraic knots are of finite type and of finite local type
(\cite{ShCa}).\newline

\noindent(2) \ \ Let $V$ be a reduced projective hypersurface of degree $d$ in $\mathbb{CP}^{n+1}$
($n \geq 1$).
Choose a Whitney stratification $\mathcal{S}$ of $V$. Recall that there is such a stratification where
strata are pure dimensional locally closed algebraic subsets with a finite number of irreducible nonsingular components.
Together with the hypersurface complement, $\mathbb{CP}^{n+1} - V$, this gives a stratification of the pair
$(\mathbb{CP}^{n+1},V)$, in which $\mathcal{S}$ is the set of singular strata.
All  links of  strata of the pair $(\mathbb{CP}^{n+1},V)$ are algebraic, hence of finite (homological) type,
so $V \subset \mathbb{CP}^{n+1}$ is of finite
local type (see \cite{ShCa}, Proposition 2.2). We choose a generic hyperplane $H$ in $\mathbb{CP}^{n+1}$,
 i.e.  transversal to all the strata of $V$, and consider the induced  stratification on the pair $(\mathbb{CP}^{n+1},V
 \cup H)$, with (open) strata of the form
 $S - S \cap H$, $S \cap H$ and $H - V \cap H$, for  $S \in \mathcal{S}$.
We call $H$ 'the hyperplane at infinity' and say that '$V$ is transversal to the hyperplane at infinity'.
Following \cite{ShCa}, we define a local system $\mathcal{L}_{H}$ on
 $\mathbb{CP}^{n+1} - (V \cup H)$,
with stalk $\Gamma := \Gamma_{\mathbb{Q}}= \mathbb{Q} [t, t^{-1}]$ and action by an element
$\alpha\in\pi_1(\mathbb{CP}^{n+1} - V \cup H)$ determined by multiplication by
$t^{\text{lk}(V \cup -dH,\alpha)}$. Here $\text{lk}(V \cup -dH,\alpha)$ is the linking number of $\alpha$ with the divisor
$V \cup -dH$ of $\mathbb{CP}^{n+1}$.
Then, (using a triangulation  of the projective space) $V \cup H$ is a
(PL) sub-pseudomanifold of $\mathbb{CP}^{n+1}$ and the intersection complex
$IC_{\bar p} ^{\bullet}:= IC_{\bar p} ^{\bullet} (\mathbb{CP}^{n+1}, \mathcal{L}_{H})$ is defined for any
(super-)perversity ${\bar p}$. The modules
$$IH_i^{\bar m}(\mathbb{CP}^{n+1};\mathcal{L}_{H}):=
\mathbb{H}^{-i}(\mathbb{CP}^{n+1};IC_{\bar m} ^{\bullet})$$
will be called the \emph{intersection Alexander modules of the hypersurface $V$}. Note that these modules are of finite
type over $\Gamma$, since $IC_{\bar p} ^{\bullet}$ is cohomologically constructible (\cite{B}, V.3.12) and $\mathbb{CP}^{n+1}$ is
compact (see \cite{B}, V.3.4.(a), V.10.13).

It will be useful to describe the links of the pair $(\mathbb{CP}^{n+1},V \cup H)$ in terms of those of
$(\mathbb{CP}^{n+1},V)$. Because of the transversality assumption, there are  stratifications $\{ Z_{i} \}$ of
$(\mathbb{CP}^{n+1},V)$ and $\{ Y_{i} \}$ of $(\mathbb{CP}^{n+1},V\cup H)$ with
$Y_{i}=Z_{i} \cup (Z_{i+2} \cap H)$ (here the indices indicate the real dimensions).
The link pair of a point $y \in (Y_{i} - Y_{i-1}) \cap H = (Z_{i+2} - Z_{i+1}) \cap H$ in
$(\mathbb{CP}^{n+1},V\cup H)$ is
$(G,F)=(S^{1} \ast G_{1},(S^{1} \ast F_{1}) \cup G_{1})$, where $(G_{1},F_{1})$ is the link pair of
$y \in Z_{i+2} - Z_{i+1}$ in $(\mathbb{CP}^{n+1},V)$.
Points in $V - V \cap H$ have the same link pairs in $(\mathbb{CP}^{n+1},V)$ and $(\mathbb{CP}^{n+1},V \cup H)$.
Finally, the link pair at any point in  $H - V \cap H$ is $(S^{1}, \emptyset)$. (For details, see \cite{ShCa}).

By Lemma 2.3.1 of \cite{ShCa}, $V \cup H \subset \mathbb{CP}^{n+1}$ is of finite local type.
Hence, by Theorem 3.3 of \cite{ShCa}, we have the following isomorphism (the \emph{superduality isomorphism}):
$$IC_{\bar m} ^{\bullet} \cong \mathcal{D} {IC_{\bar l} ^{\bullet}}^{op} [2n+2]$$
(here $A^{op}$ is the $\Gamma$-module obtained from the $\Gamma$-module $A$ by composing all module structures with the
involution $t \to t^{-1}$.)
Recall that the middle and logarithmic perversities are defined as: ${\bar m}(s) = [(s-1)/2]$ and
${\bar l}(s) = [(s+1)/2]$. Note ${\bar m}(s) + {\bar l}(s) =s-1$, i.e.  ${\bar m}$ and ${\bar l}$ are superdual perversities.\newline

\noindent(3) \ \ With the notations from $\S 2$, we have an  isomorphism of $\Gamma$-modules:
$$H_i(\mathcal{U};\mathcal{L}_{H}) \cong H_i(\mathcal{U}^c;\mathbb{Q})$$
where $\mathcal{L}_{H}$ is, as above, the local coefficient system on $\mathcal{U}=\mathbb{CP}^{n+1} - (V \cup H)$
defined by the representation
$\mu : \pi_1 (\mathcal{U}) \to \text{Aut}(\Gamma)=\Gamma^*$, $\mu(\alpha)=t^{\text{lk}(\alpha,V \cup -dH)}$.

Indeed, $\mathcal{U}^c$ is  the covering associated to the kernel of the linking number homomorphism
$lk : \pi_1(\mathcal{U}) \to \mathbb{Z}$, $\alpha \mapsto \text{lk}(\alpha,V \cup -dH)$, and note that $\mu$
factors through $lk$, i.e. $\mu$ is the composition
$ \pi_1(\mathcal{U}) \overset{lk}{\to} \mathbb{Z} \to \Gamma^*$, with the second homomorphism mapping $1$ to $t$.
 Thus $\text{Ker}(lk) \subset \text{Ker}(\mu)$. By definition, $H_*(\mathcal{U};\mathcal{L}_{H})$ is the
 homology of the chain complex $C_*(\mathcal{U};\mathcal{L}_H)$ defined by the equivariant
 tensor product: $C_*(\mathcal{U};\mathcal{L}_H):=C_*(\mathcal{U}^c) \otimes_{\mathbb{Z}} \Gamma$,
 where $\mathbb{Z}$ stands for the group of covering transformations of $\mathcal{U}^c$ (see \cite{Di2}, page 50).
 Since $\Gamma=\mathbb{Q}[\mathbb{Z}]$,  the chain complex
$C_*(\mathcal{U}^c) \otimes_{\mathbb{Z}} \Gamma$ is clearly isomorphic to the complex
$C_*(\mathcal{U}^c) \otimes \mathbb{Q}$, and the claimed isomorphism follows.
For a similar argument, see also \cite{Ha}, Example 3H.2.\newline

\noindent(4) \ \ This is also a convenient place to point out the following fact: because ${\bar m}(2)=0$,
the allowable zero- and one-chains (\cite{GM1}) are
those which lie in  $\mathbb{CP}^{n+1} - (V \cup H)$. Therefore,
$$IH_0^{\bar m}(\mathbb{CP}^{n+1};\mathcal{L}_{H}) \cong H_0(\mathbb{CP}^{n+1} - (V \cup H);\mathcal{L}_{H})=\Gamma /(t-1),$$
where the second isomorphism follows from the identification of the homology of $\mathbb{CP}^{n+1} - (V \cup H)$
with local coefficient system $\mathcal{L}_H$ with the rational homology (viewed as a $\Gamma$-module) of the infinite cyclic
 cover  $\mathcal{U}^c$  of $\mathbb{CP}^{n+1} - (V \cup H)$, defined by the linking number homomorphism.\newline


\subsection{Relation with the classical Alexander modules of the complement}
Let $V$ be a degree $d$, reduced, n-dimensional projective hypersurface, which is
transversal to the hyperplane at infinity, $H$. We are aiming to show that, in our setting,
\emph{the intersection Alexander modules of a hypersurface coincide with the classical Alexander
modules of the hypersurface complement}.
 The key fact will be the following characterization of the support of
the intersection homology complex $IC_{\bar m} ^{\bullet}$:

\begin{lemma} There is a quasi-isomorphism:
 $$IC_{\bar m} ^{\bullet} (\mathbb{CP}^{n+1} , \mathcal{L}_{H})_{|{V \cup H}} \cong 0$$
\end{lemma}

\begin{proof}
It suffices to show the vanishing of stalks of the complex $IC_{\bar m} ^{\bullet}$ at points in strata of $V \cup H$.
We will do this in two steps:\newline
\noindent\textbf{Step 1.}  $$IC_{\bar m} ^{\bullet} (\mathbb{CP}^{n+1} , \mathcal{L}_{H})_{|{H}} \cong 0$$

The link pair of $H - V \cap H$ is $(S^{1},\emptyset)$ and this maps to $t^{-d}$ under $\mathcal{L}_{H}$,
therefore the stalk of
$IC_{\bar m} ^{\bullet}(\mathbb{CP}^{n+1} , \mathcal{L}_{H})$ at a point in this stratum is zero. Indeed (cf. \cite{B},
V.3.15), for $x \in H - V \cap H$:
\begin{equation*}
\mathcal{H} ^{q} (IC_{\bar m} ^{\bullet}(\mathbb{CP}^{n+1} , \mathcal{L}_{H}))_{x} \cong
\begin{cases}
0 \ , & q > -2n-2\\
IH^{\bar m} _{-q-(2n+1)} (S^{1} ; \mathcal{L}), & q \leq -2n-2,
\end{cases}
\end{equation*}
and note that $IH^{\bar m} _{j} (S^{1} ; \mathcal{L}) \cong 0$ unless $j=0$.

Next, consider the link pair $(G,F)$ of a point $x \in S \cap H$, $S \in \mathcal{S}$.
Let the real codimension of $S$ in $\mathbb{CP}^{n+1}$
be $2k$. Then the codimension of $S \cap H$ is $2k+2$ and $\text{dim}(G)=2k+1$.
The stalk at $x \in S \cap H$ of the intersection
homology complex $IC_{\bar m} ^{\bullet}(\mathbb{CP}^{n+1} , \mathcal{L}_{H})$ is given by the local calculation
formula (\cite{B}, (3.15)):
\begin{equation*}
\mathcal{H} ^{q} (IC_{\bar m} ^{\bullet}(\mathbb{CP}^{n+1} , \mathcal{L}_{H}))_{x} \cong
\begin{cases}
0 \ , & q > k-2n-2\\
IH^{\bar m} _{-q-(2n-2k+1)} (G ; \mathcal{L}), & q \leq k-2n-2.
\end{cases}
\end{equation*}

We claim that :
$$IH_{i} ^{\bar m}(G;\mathcal{L}) = 0 \ , \ i \geq k+1$$
Then, by setting $i=-q-(2n-2k+1)$ , we obtain that
$IH^{\bar m} _{-q-(2n-2k+1)} (G ; \mathcal{L})=0 \ \text{for} \ q \leq k-2n-2$, and
therefore:
$$\mathcal{H} ^{q} (IC_{\bar m} ^{\bullet}(\mathbb{CP}^{n+1} , \mathcal{L}_{H}))_{x}=0$$

In order to prove the claim, we use arguments similar to those used in \cite{ShCa}, p. 359-361.
Recall that  $(G,F)$ is of the form $(S^{1} \ast G_{1},(S^{1} \ast F_{1}) \cup G_{1})$,
where $(G_{1},F_{1})$ is the link pair of $x \in S$ in $(\mathbb{CP}^{n+1},V)$ (or equivalently, the link of
$S \cap H$ in $(H,V \cap H) $).
The restriction of $\mathcal{L}_H$ to $G - F$ is given by sending $\alpha \in \pi_{1}(G - F)$ to
$t^{\ell(S^{1} \ast F_{1} - d G_{1},\alpha)} \in Aut(\Gamma)$.
Let $\mathcal{IC}=IC_{\bar m} ^{\bullet}(G,\mathcal{L})$. The link of the codimension two stratum $G_{1} - G_{1} \cap
(S^{1} \ast F_{1})$ of $G$ is a circle that maps to $t^{-d}$ under $\mathcal{L}$; hence by the
stalk cohomology formula (\cite{B}, (3.15)), for $y \in G_{1} - G_{1} \cap (S^{1} \ast F_{1})$, we have:
\begin{equation*}
\mathcal{H} ^{i} (\mathcal{IC})_{y} \cong
\begin{cases}
IH^{\bar m} _{-i-(\text{dim}G_1 +1)} (S^{1} ; \Gamma), & i \leq -\text{dim}G\\
0 \ , & i > {\bar m}(2)-\text{dim}G=-\text{dim}G.
\end{cases}
\end{equation*}
Since $IH^{\bar m} _{-i-(\text{dim}G_{1}+1)} (S^{1} ; \Gamma) \cong H_{-i-(\text{dim}G-1)} (S^{1} ; \Gamma)=0$ for
$-i-(\text{dim}G-1) \neq 0$, i.e., for $i \neq 1-\text{dim}G$, we obtain that:
$$\mathcal{IC}|_{G_{1} - G_{1} \cap (S^{1} \ast F_{1})} \cong 0$$
Moreover, $G_{1}$ is a locally flat submanifold of $G$ and
intersects $S^{1} \ast F_{1}$ transversally. Hence the link pair
in $(G,F)$ of a stratum of $G_{1} \cap (S^{1} \ast F_{1})$ and the
restriction of $\mathcal{L}$ will have the same form as links of
strata of $V \cap H$ in $(\mathbb{CP}^{n+1},V \cup H)$. Thus, by
induction on dimension we obtain:
$$\mathcal{IC}|_{G_{1} \cap (S^{1} \ast F_{1})} \cong 0$$
Therefore, $\mathcal{IC}|_{G_{1}} \cong 0$.
Thus, denoting by $i$ and $j$ the inclusions of $G - G_{1}$ and $G_{1}$, respectively,
the distinguished triangle: $$Ri_{!}i^{\ast}\mathcal{IC} \to \mathcal{IC} \to Rj_{\ast}\mathcal{IC}|_{G_{1}} \to$$
upon applying the compactly supported hypercohomology functor, yields the isomorphisms:
$${^{c}IH_{i} ^{\bar m}} (G - G_{1};\mathcal{L}) \cong IH_{i} ^{\bar m} (G;\mathcal{L})$$
We have: $$(G - G_{1}, F - F \cap G_{1}) \cong
(c^{\circ}G_{1} \times S^{1} , c^{\circ}F_{1} \times S^{1})$$
and $\mathcal{L}$ is given on
$$(c^{\circ}G_{1} - c^{\circ}F_{1}) \times S^{1} \cong (G_{1} - F_{1}) \times \mathbb{R}
\times S^{1}$$
by sending $\alpha \in \pi_{1}(G_{1} - F_{1})$ to the multiplication by $t^{\ell(F_{1},\alpha)}$, and a generator
of $\pi_{1}(S^{1})$ to $t^{-d}$.

We denote by $\mathcal{L}_{1}$ and $\mathcal{L}_{2}$ the restrictions of $\mathcal{L}$ to $c^{\circ}G_{1}$ and $S^{1}$
respectively.
Note that $IH_{b} ^{\bar m} (S^{1};\mathcal{L}_{2})=0$ unless $b=0$, in which case it is isomorphic to $\Gamma/t^{d}-1$.
Therefore, by the Kunneth formula (\cite{GM4}), we have:
\[
\begin{aligned}
{^{c}IH_{i} ^{\bar m}} (c^{\circ}G_{1} \times S^{1};\mathcal{L}) &\cong
\{{^{c}IH_{i} ^{\bar m}} (c^{\circ}G_{1};\mathcal{L}_{1}) \otimes IH_{0} ^{\bar
m} (S^{1};\mathcal{L}_{2})\}\\
  &\oplus
\{{^{c}IH_{i-1} ^{\bar m}} (c^{\circ}G_{1};\mathcal{L}_{1}) \ast IH_{0} ^{\bar m} (S^{1};\mathcal{L}_{2})\}
\end{aligned}
\]

Lastly, the formula for the compactly supported intersection homology of a cone yields (\cite{B}, \cite{Ki}, \cite{F}):
${^{c}IH_{i} ^{\bar m}} (c^{\circ}G_{1};\mathcal{L}_{1})=0$ for $i \geq \text{dim}G_{1} - {\bar m} (\text{dim}G_{1}+1)=2k-1-{\bar
m}(2k)=k$ (as well as  ${^{c}IH_{i} ^{\bar m}} (c^{\circ}G_{1};\mathcal{L}_{1})=IH_{i} ^{\bar m} (G_{1};\mathcal{L}_{1})$,
for $i < \text{dim}G_{1} - {\bar m} (\text{dim}G_{1}+1)=k$)
and, consequently: ${^{c}IH_{i-1} ^{\bar m}} (c^{\circ}G_{1};\mathcal{L}_{1})=0$ for $i \geq k+1$.

Altogether,
$$IH_{i} ^{\bar m} (G;\mathcal{L}) \cong {^{c}IH_{i} ^{\bar m}} (G - G_{1};\mathcal{L})
\cong {^{c}IH_{i} ^{\bar m}} (c^{\circ}G_{1} \times S^{1};\mathcal{L})=0 \ \text{for} \ i \geq k+1$$
as claimed.\newline

\noindent\textbf{Step 2.} $$IC_{\bar m} ^{\bullet} (\mathbb{CP}^{n+1} , \mathcal{L}_{H})_{|{V}} \cong 0$$

It suffices to show the vanishing of stalks of the complex $IC_{\bar m} ^{\bullet}$ at points in strata of the form $S -
S \cap H$ of the affine part, $V_{aff}$, of $V$. Note that, assuming $S$ connected, the link pair of $S -
S \cap H$ in $(\mathbb{CP}^{n+1},V \cup H)$ is the same as its link pair in $(\mathbb{C}^{n+1},V_{aff})$ with the induced
stratification, or the link pair
of $S$ in $(\mathbb{CP}^{n+1},V)$.
Let $x \in S-S\cap H$ be a point in an affine stratum of complex dimension $s$. The stalk cohomology calculation yields:
\begin{equation*}
\mathcal{H} ^{q} (IC_{\bar m} ^{\bullet})_{x} \cong
\begin{cases}
IH^{\bar m} _{-q-(2s+1)} (S_x ^{2n-2s+1} ; \Gamma), & q \leq -n-s-2\\
0 \ , & q > -n-s-2,
\end{cases}
\end{equation*}
where $(S_x ^{2n-2s+1}, K_x)$ is the link pair of the component containing $x$.

To obtain the desired vanishing, it suffices to prove that:
$IH^{\bar m} _{j} (S_x ^{2n-2s+1} ; \Gamma) \cong 0$ for $j \geq n-s+1$. We will show this in
the following:

\begin{lemma} If $S$ is an $s$-dimensional stratum of $V_{aff}$ and $x$ is a point in $S$, then the intersection homology groups of
its link pair $(S_x ^{2n-2s+1}, K_x)$ in $(\mathbb{C}^{n+1},V_{aff})$ are characterized by the following properties:
$$IH^{\bar m} _{j} (S_x ^{2n-2s+1} ; \Gamma) \cong 0 \ \ , \ \ j \geq n-s+1$$
$$IH^{\bar m} _{j} (S_x ^{2n-2s+1} ; \Gamma) \cong H_j (S_x ^{2n-2s+1} - K_x; \Gamma) \ \ , \ \ j \leq n-s.$$
\end{lemma}

\noindent(here $\Gamma$ denotes the local coefficient system on the knot complement, with stalk $\Gamma$ and action of an element
$\alpha$ in the fundamental group of the complement given by multiplication by $t^{\text{lk}(\alpha,K)}$; \cite{ShCa},
\cite{F})

\emph{Note.} The same property holds for link pairs of strata $S \in \mathcal{S}$ of a stratification of the pair
$(\mathbb{CP}^{n+1},V)$ since all of these are algebraic knots and have associated Milnor fibrations.

\emph{Proof of lemma.} We will prove the above claim by induction down on the dimension of singular strata of the pair
$(\mathbb{C}^{n+1},V_{aff})$. To start the induction, note that the link pair of a component of the dense open subspace of $V_{aff}$
(i.e. for $s=n$) is
a circle  $(S^{1}, \emptyset)$, that maps to $t$ under $\mathcal{L}_{H}$. Moreover, the (intersection) homology groups
$IH_{i}^{\bar m}(S^{1};\Gamma) \cong H_{i}(S^{1};\Gamma)$ are zero, except for $i=0$, hence the claim is trivially satisfied in this case.

Let $S$ be an $s$-dimensional stratum of $V_{aff}$ and let $x$ be a point in $S$. Its link pair
$(S_x ^{2n-2s+1}, K_x)$ in $(\mathbb{C}^{n+1},V_{aff})$ is a singular algebraic knot,
with a topological stratification induced by that of
$(\mathbb{C}^{n+1},V_{aff})$. The link pairs of strata of $(S_x ^{2n-2s+1}, K_x)$ are also link pairs of higher dimensional
strata of $(\mathbb{C}^{n+1},V_{aff})$ (see for example \cite{F}).
Therefore, by the induction hypothesis, the claim holds for such link pairs.

Let $\mathcal{IC}=IC_{\bar m} ^{\bullet}(S_x ^{2n-2s+1},\Gamma)$ be the middle-perversity intersection cohomology complex
associated to the link pair of $S$ at $x$. In order to prove the claim, it suffices to show that its restriction to $K$ is
quasi-isomorphic to the zero complex, i.e.
$\mathcal{IC}_{|K} \cong 0$. Then the lemma will follow from the long exact sequence of compactly supported hypercohomology
and from the fact that the fiber $F_x$ of the Milnor fibration associated with the algebraic knot
$(S_x ^{2n-2s+1}, K_x)$ has the homotopy
type of an $(n-s)$-dimensional complex (\cite{Mi}, Theorem 5.1) and is homotopy equivalent to the infinite cyclic  covering $\widetilde{S_x ^{2n-2s+1} - K_x}$
of the knot complement, defined by the linking number homomorphism. More precisely, we obtain the isomorphisms:
$$IH^{\bar m} _{j} (S_x ^{2n-2s+1} ; \Gamma) \cong H_j (S_x ^{2n-2s+1} - K_x; \Gamma) \cong
H_j (\widetilde{S_x ^{2n-2s+1} - K_x}; \mathbb{Q}) \cong H_j (F_x;\mathbb{Q}).$$

Let $K' \supset K''$ be two consecutive terms in the filtration of $(S_x ^{2n-2s+1}, K_x)$. Say
$\text{dim}_{\mathbb{R}}(K')=2n-2r-1$, $r \geq s$. The stalk of $\mathcal{IC}$ at a point $y \in K'-K''$ is given by the
following formula:
\begin{equation*}
\mathcal{H} ^{q} (\mathcal{IC})_{y} \cong
\begin{cases}
IH^{\bar m} _{-q-(2n-2r)} (S_y ^{2r-2s+1} ; \Gamma), & q \leq -2n+s+r-1\\
0 \ , & q > -2n+s+r-1,
\end{cases}
\end{equation*}
where $(S_y ^{2r-2s+1}, K_y)$ is the link pair in $(S_x ^{2n-2s+1},K_x)$ of the component of $K'-K''$ containing $y$.
Since $(S_y ^{2r-2s+1}, K_y)$ is also the link pair of a higher dimensional
stratum of $(\mathbb{C}^{n+1},V_{aff})$,
 the induction hypothesis yields: $IH^{\bar m} _{-q-(2n-2r)} (S_y ^{2r-2s+1} ; \Gamma) \cong 0 \ \text{if} \  q \leq -2n+s+r-1$.
\end{proof}

\begin{remark}
The proof of Step 1 the previous lemma  provides a way of computing the modules  $IH_{i} ^{\bar m} (G;\mathcal{L})$, $i \leq k$, for
$G \cong S^{2k+1} \cong S^{1}*G_1$ the link of an $(n-k)$-dimensional stratum $S \cap H$, $S \in \mathcal{S}$:
$$IH_{k} ^{\bar m} (G;\mathcal{L}) \cong
{^{c}IH_{k-1} ^{\bar m}} (c^{\circ}G_{1};\mathcal{L}_{1}) \ast IH_{0} ^{\bar m} (S^{1};\mathcal{L}_{2})
\cong IH_{k-1} ^{\bar m} (G_{1};\mathcal{L}_{1}) \ast IH_{0} ^{\bar m} (S^{1};\mathcal{L}_{2})$$
and, for $i <k$:
$$IH_{i} ^{\bar m} (G;\mathcal{L}) \cong
\{ IH_{i} ^{\bar m} (G_{1};\mathcal{L}_{1}) \otimes IH_{0} ^{\bar m} (S^{1};\mathcal{L}_{2}) \} \oplus
\{ IH_{i-1} ^{\bar m} (G_{1};\mathcal{L}_{1}) \ast IH_{0} ^{\bar m} (S^{1};\mathcal{L}_{2}) \}.$$
The above formulas, as well as the claim of the first step of the previous lemma, can also be obtained from the formula for the
intersection homology of a join (\cite{GM4}, Proposition 3), applied to  $G \cong G_1 \ast S^{1}$.

If we denote by $I\gamma ^{\bar m}_{i}(G):=\text{order} \ IH_{i} ^{\bar m} (G;\mathcal{L})$,
the intersection Alexander polynomial of the link pair $(G,F)$ (see \cite{F}), then we obtain:
$$I\gamma ^{\bar m}_{k}(G)=\text{gcd} (I\gamma ^{\bar m}_{k-1}(G_1) , t^{d}-1)$$
$$I\gamma ^{\bar m}_{i}(G)=\text{gcd} (I\gamma ^{\bar m}_{i}(G_1) , t^{d}-1) \times
\text{gcd} (I\gamma ^{\bar m}_{i-1}(G_1) , t^{d}-1) \ , \ i <k \ .$$
In particular, since $I\gamma ^{\bar m}_{0}(G_1) \sim t-1$ (\cite{F}, Corollary 5.3), we have: $I\gamma ^{\bar m}_{0}(G) \sim
t-1$.

Note that the superduality isomorphism (\cite{ShCa}, Corollary 3.4) yields the isomorphism: $IH_{j} ^{\bar l} (G;\mathcal{L}) \cong
IH_{2k-j} ^{\bar m} (G;\mathcal{L})^{op}$. Hence $IH_{j} ^{\bar l} (G;\mathcal{L}) \cong 0$ if $j<k$.

>From the above considerations, the zeros of the polynomials $I\gamma ^{\bar m}_{i}(G)$ and $I\gamma ^{\bar l}_{i}(G)$
(in the non-trivial range) are all roots of unity of order $d$.
\end{remark}

\begin{cor}
If $V$ is an $n$-dimensional  reduced  projective hypersurface, transversal to the hyperplane at infinity,
then the intersection Alexander modules of $V$ are isomorphic to the classical Alexander modules of the hypersurface
complement, i.e.
$$IH_{*} ^{\bar m} (\mathbb{CP}^{n+1};\mathcal{L}_{H}) \cong H_* (\mathbb{CP}^{n+1} - V \cup H; \mathcal{L}_H) \cong
H_{*}(\mathcal{U} ^{c};\mathbb{Q})$$
\end{cor}
\begin{proof}
The previous lemma and the hypercohomology spectral sequence yield:
$$\mathbb{H} ^{i}  (V \cup H;IC_{\bar m} ^{\bullet}) \cong 0$$
Let $u$ and $v$ be the inclusions of $\mathbb{CP}^{n+1} - (V \cup H)$ and resp. $V \cup H$ into $\mathbb{CP}^{n+1}$. The
distinguished triangle $u_!u^* \to id \to v_*v^* \overset{[1]}{\to}$, upon applying the hypercohomology
functor, yields the long exact sequence:
\[
\begin{aligned}
\cdots \to \mathbb{H}_c ^{-i-1}(V \cup H;IC_{\bar m} ^{\bullet}) \to
\mathbb{H}^{-i} _{c} (\mathbb{CP}^{n+1} - (V \cup H);IC_{\bar m} ^{\bullet}) \to
\\
\to \mathbb{H}_c ^{-i} (\mathbb{CP}^{n+1};IC_{\bar m} ^{\bullet})
\to \mathbb{H}_c ^{-i}(V \cup H;IC_{\bar m} ^{\bullet}) \to \cdots
\end{aligned}
\]
Therefore, we obtain the isomorphisms:
\[
\begin{aligned}
IH_{i} ^{\bar m} (\mathbb{CP}^{n+1};\mathcal{L}_{H})&:=\mathbb{H}^{-i} (\mathbb{CP}^{n+1};IC_{\bar m} ^{\bullet})\\
&\cong \mathbb{H}^{-i} _{c} (\mathbb{CP}^{n+1} - (V \cup H);IC_{\bar m} ^{\bullet})\\
&= {^{c}IH_{i} ^{\bar m}}  (\mathbb{CP}^{n+1} - (V \cup H);\mathcal{L}_{H})\\
&\cong H_{i} (\mathbb{CP}^{n+1} - (V \cup H);\mathcal{L}_H)\\
&\cong H_{i}(\mathcal{U} ^{c};\mathbb{Q})
\end{aligned}
\]
\end{proof}

 Our next goal is to show that, in our settings, the Alexander modules of the hypersurface complement,
 $H_i (\mathcal{U}^c; \mathbb{Q})$, are torsion $\Gamma$-modules if $i \leq n$. Based on the above corollary, it suffices to
 show this for the modules $IH_{i} ^{\bar m} (\mathbb{CP}^{n+1};\mathcal{L}_{H})$, $i \leq n$.

 We will need the following:
\begin{lemma}
$$IH_{i} ^{\bar l} (\mathbb{CP}^{n+1};\mathcal{L}_{H}) \cong 0 \ \ \text{for} \ \ i \leq n$$
\end{lemma}

\begin{proof}
Let $u$ and $v$ be the inclusions of $\mathbb{CP}^{n+1} - V \cup H$
and respectively $V \cup H$ into
$\mathbb{CP}^{n+1}$. Since $v^{*}IC_{\bar m} ^{\bullet} \cong 0$, by superduality we obtain:
$0 \cong v^{*}\mathcal{D}{IC_{\bar l} ^{\bullet} [2n+2]}^{op} \cong
\mathcal{D}v^{!}{IC_{\bar l} ^{\bullet} [2n+2]}^{op}$, so
$v^{!}IC_{\bar l} ^{\bullet} \cong 0$. Hence the
distinguished triangle:
$$v_{*}v^{!}IC_{\bar l} ^{\bullet} \to IC_{\bar l} ^{\bullet} \to u_{*}u^{*}IC_{\bar l} ^{\bullet} \overset{[1]}{\to}$$
upon applying the hypercohomology functor, yields the isomorphism:
$$IH_{i} ^{\bar l} (\mathbb{CP}^{n+1};\mathcal{L}_{H}) \cong IH_{i} ^{\bar l} (\mathbb{CP}^{n+1} - V \cup H;\mathcal{L}_H)
\cong H^{BM} _i (\mathbb{CP}^{n+1} - V \cup H;\mathcal{L}_H)$$
where $H^{BM} _*$ denotes the Borel-Moore homology.
By Artin's vanishing theorem (\cite{S}, Example 6.0.6), the latter module is $0$ for $i < n+1$,
since $\mathbb{CP}^{n+1} - V \cup H$ is a Stein space of dimension $n+1$.
\end{proof}

\begin{remark}
Recall that the \emph{peripheral complex} $\mathcal{R}^{\bullet}$, associated to the finite local type
embedding $V \cup H \subset \mathbb{CP}^{n+1}$,
is defined by the distinguished triangle (\cite{ShCa}):
$$IC_{\bar m} ^{\bullet} \to IC_{\bar l} ^{\bullet}  \to \mathcal{R}^{\bullet} \overset{[1]}{\to}$$
Moreover, $\mathcal{R}^{\bullet}$ is a perverse (in the sense considered in \cite{ShCa}),
self-dual (i.e., $\mathcal{R}^{\bullet} \cong \mathcal{D} {\mathcal{R}^{\bullet}}^{op}[2n+3]$),
torsion sheaf on $\mathbb{CP}^{n+1}$ (i.e., the stalks of its cohomology sheaves
are torsion modules). All these properties are preserved by restriction to open sets.
\end{remark}

By applying the hypercohomology functor to the triangle defining the peripheral complex $\mathcal{R}^{\bullet}$, and using
the vanishing of the previous lemma, we obtain the following:

\begin{prop}
The natural maps:
$$\mathbb{H}^{-i-1} (\mathbb{CP}^{n+1};\mathcal{R}^{\bullet}) \to IH_{i} ^{\bar m} (\mathbb{CP}^{n+1};\mathcal{L}_{H})$$
are isomorphisms for all $i \leq n-1$ and epimorphism for $i=n$.
\end{prop}

Now, since $\mathcal{R}^{\bullet}$ is a torsion sheaf (having finite dimensional rational vector spaces as stalks),
the spectral sequence for
hypercohomology implies that the groups $\mathbb{H}^{q} (\mathbb{CP}^{n+1};\mathcal{R}^{\bullet})$, $q \in \mathbb{Z}$, are also finite
dimensional rational vector spaces, thus torsion $\Gamma$-modules. Therefore, the above proposition
yields the following:
\begin{cor}
Let $V \subset \mathbb{CP}^{n+1}$ be a reduced, $n$-dimensional projective hypersurface,
transversal to the
hyperplane at infinity.
Then for any $i \leq n$, the module $IH_{i} ^{\bar m} (\mathbb{CP}^{n+1};\mathcal{L}_{H}) \cong
H_i (\mathcal{U}^c;\mathbb{Q})$ is a finitely generated torsion $\Gamma$-module.
\end{cor}

\emph{Note.} (1) For $i \leq n$, $H_i (\mathcal{U}^c;\mathbb{Q})$ is actually a finite dimensional rational vector space, thus
its order coincides with the characteristic polynomial of the $\mathbb{Q}$-linear map induced by a generator of the group of covering
transformations (see \cite{Mi2}).\newline
(2) Lemma 1.5 of \cite{Li} asserts that if $k$ is the dimension of the singular locus of $V$
and $n-k \geq 2$, then
$H_{i}(\mathcal{U} ^{c};\mathbb{Q}) \cong 0$ for $1 \leq i < n-k$ (here we use the fact that, if $n-k \geq 2$,
 the infinite cyclic cover of $\mathbb{CP}^{n+1} - V \cup H$ is the universal cover).

\begin{defn}
For $i \leq n$, we denote by $\delta_i (t)$ the polynomial associated to the torsion module $H_i
(\mathcal{U}^c;\mathbb{Q})$, and call it the $i$-th global Alexander polynomial of the hypersurface $V$.
These polynomials will be well-defined up to multiplication by $ct^k$, $c \in \mathbb{Q}$.
\end{defn}

As a consequence of the previous corollary, we may calculate the rank of the free $\Gamma$-module
$H_{n+1}(\mathcal{U}^c;\mathbb{Q})$ in terms of the Euler characteristic of the complement:
\begin{cor}
Let $V \subset \mathbb{CP}^{n+1}$ be a reduced, $n$-dimensional projective hypersurface, in general position at infinity.
Then the $\Gamma$-rank of $H_{n+1}(\mathcal{U}^c;\mathbb{Q})$ is expressed in terms of the
Euler characteristic $\chi (\mathcal{U})$ of the complement by the formula:
$$(-1)^{n+1} \chi(\mathcal{U})= \text{rank}_{\Gamma} H_{n+1}(\mathcal{U}^c;\mathbb{Q})$$
\end{cor}
\begin{proof}
The equality follows from the above corollary, from the fact that for $q>n+1$ the Alexander modules
$H_{q}(\mathcal{U}^c;\mathbb{Q})$ vanish, and from the formula 2.10(v) of \cite{Di3}:
$$\chi (\mathcal{U})=\sum_q (-1)^{q} \text{rank}_{\Gamma} H_{q}(\mathcal{U}^c;\mathbb{Q})$$
\end{proof}


\section{The Main Theorems}
We will now state and prove the main theorems of this paper. These
results are generalizations of the ones obtained by A. Libgober
(\cite{Li}, \cite{Li3}, \cite{Li4}) in the case of hypersurfaces
with isolated singularities, and will lead to results on the
monodromy of the Milnor fiber of a projective hypersurface
arrangement, similar to those obtained by Libgober (\cite{Li100}),
Dimca (\cite{Di2}, \cite{COD}) etc (see \S 5).

The first theorem provides a characterization of the zeros of
global Alexander polynomials. For hypersurfaces with only isolated
singularities, it specializes to Corollary 4.8 of \cite{Li}. It
also gives a first obstruction on the prime divisors of the global
Alexander polynomials of hypersurfaces:

\begin{thm}
If $V$ is an n-dimensional reduced projective hypersurface of
degree $d$, transversal to the hyperplane at infinity, then for $i
\leq n$, any root of the  global Alexander polynomial
$\delta_{i}(t)$ is a root of unity of order $d$.
\end{thm}
\begin{proof}
Let $k$ and $l$ be the inclusions of $\mathbb{C}^{n+1}$ and respectively $H$ into $\mathbb{CP}^{n+1}$. For a fixed
perversity $\bar p$, we will denote the
intersection complexes $IC_{\bar p} ^{\bullet} (\mathbb{CP}^{n+1} , \mathcal{L}_H)$  by $IC_{\bar p} ^{\bullet}$.
We will also drop the letter $R$ when using right derived functors.
The distinguished
triangle:
$l_* l^! \to id \to k_* k^* \overset{[1]}{\to}$,
upon applying the hypercohomology functor, yields the
following exact sequence:
\[
\begin{aligned}
\cdots \to \mathbb{H}^{-i}_{H}(\mathbb{CP}^{n+1};IC_{\bar m} ^{\bullet}) \to
IH_{i} ^{\bar m} (\mathbb{CP}^{n+1};\mathcal{L}_{H}) \to \mathbb{H}^{-i}(\mathbb{C}^{n+1};k^*IC_{\bar m} ^{\bullet}) \to \\
\to
\mathbb{H}^{-i+1}_{H}(\mathbb{CP}^{n+1};IC_{\bar m} ^{\bullet}) \to \cdots
\end{aligned}
\]
Note that the complex $k^*IC_{\bar m} ^{\bullet}[-n-1]$
is perverse with respect to the middle perversity (since $k$ is the open
inclusion and the functor $k^*$ is t-exact; \cite{BBD}). Therefore, by Artin's
vanishing theorem for perverse sheaves (\cite{S}, Corollary 6.0.4), we obtain:
$$\mathbb{H}^{-i} (\mathbb{C}^{n+1};k^*IC_{\bar m} ^{\bullet}) \cong 0 \ \ \text{for} \ \ i<n+1.$$
Hence:
$$\mathbb{H}^{-i}_{H}(\mathbb{CP}^{n+1};IC_{\bar m} ^{\bullet}) \cong
IH_{i} ^{\bar m} (\mathbb{CP}^{n+1};\mathcal{L}_{H}) \ \ \text{for} \ \ i<n, $$
and $IH_{n} ^{\bar m} (\mathbb{CP}^{n+1};\mathcal{L}_{H})$ is a quotient of
$\mathbb{H}^{-n}_{H}(\mathbb{CP}^{n+1};IC_{\bar m} ^{\bullet})$.

The superduality isomorphism $IC_{\bar m} ^{\bullet} \cong \mathcal{D} {IC_{\bar l} ^{\bullet}}^{op} [2n+2]$,
and the fact that the stalks over $H$ of the complex $IC_{\bar l} ^{\bullet}$ are torsion
$\Gamma$-modules (recall that $l^* {IC_{\bar l} ^{\bullet}} \cong l^* {\mathcal{R}^{\bullet}}$, and $\mathcal{R}^{\bullet}$
is a torsion sheaf by \cite{ShCa}), yield the isomorphisms:
\[
\begin{aligned}
\mathbb{H}^{-i} _{H}(\mathbb{CP}^{n+1};IC_{\bar m} ^{\bullet})
&= \mathbb{H}^{-i} (H;l^!IC_{\bar m} ^{\bullet}) \\
&\cong \mathbb{H}^{-i+2n+2} (H;\mathcal{D}l^* {IC_{\bar l} ^{\bullet}}^{op}) \\
&\cong Hom(\mathbb{H}^{i-2n-2} (H;l^* {IC_{\bar l} ^{\bullet}}^{op});\Gamma) \oplus
Ext(\mathbb{H}^{i-2n-1}(H;l^* {IC_{\bar l} ^{\bullet}}^{op});\Gamma)\\
&\cong Ext(\mathbb{H}^{i-2n-1} (H;l^* {IC_{\bar l} ^{\bullet}}^{op});\Gamma)\\
&\cong \mathbb{H}^{i-2n-1} (H;l^* {IC_{\bar l} ^{\bullet}}^{op})\\
&\cong \mathbb{H}^{i-2n-1} (H;l^* {\mathcal{R}^{\bullet}}^{op}).\\
\end{aligned}
\]

Then, in order to finish the proof of the theorem, it suffices to
study the order of the module $\mathbb{H}^{i-2n-1} (H;l^*
{\mathcal{R}^{\bullet}}^{op})$, for $i \leq n$, and to show that
the zeros of its associated polynomial are roots of unity of order
$d$. This  follows by using the hypercohomology spectral sequence,
since the stalks of ${\mathcal{R}^{\bullet}}^{op}$ at points of
$H$ are torsion modules whose associated polynomials have the
desired property: their zeros are roots of unity of order $d$ (see
Remark 3.3 concerning the local intersection Alexander polynomials
associated to link pairs of strata of $V \cap H$).
\end{proof}

\emph{Note.} The above theorem is also a generalization of the
following special case. If $V$ is a projective cone on a degree
$d$ reduced hypersurface $Y=\{f=0\} \subset \mathbb{CP}^n$, then
there is a $\Gamma$-module isomorphism: $H_i (\mathcal{U}^c;\mathbb{Q}) \cong
H_i(F;\mathbb{Q})$, where $F=f^{-1}(1)$ is the fiber of the global Milnor
fibration $\mathbb{C}^{n+1}-f^{-1}(0) \overset{f}{\to} \mathbb{C}^*$
associated to the homogeneous polynomial $f$, and the module structure on $H_i(F;\mathbb{Q})$
is induced by the monodromy action (see \cite{Di}, p. 106-107). Therefore the zeros of the
global Alexander polynomials of $V$ coincide with the eigenvalues
of the monodromy operators acting on the homology of $F$. Since
the monodromy homeomorphism has finite order $d$, all these eigenvalues are
roots of unity of order $d$.\newline

Next we show that the zeros of the polynomials $\delta_{i}(t)$ ($i \leq n$)
are controlled by the local data, i.e. by the local Alexander polynomials of link pairs
associated to singular strata  contained in some fixed component of $V$, in
a stratification of the pair $(\mathbb{CP}^{n+1},V)$.
This is an extension to the case of general
singularities of a result due to A. Libgober (\cite{Li}, Theorem
4.3; \cite{Li4}, Theorem 4.1.a), which gives a similar fact for
the hypersurfaces with isolated singularities.

\begin{thm}
Let $V$ be a reduced hypersurface in $\mathbb{CP}^{n+1}$, which is
transversal to the hyperplane at infinity, $H$. Fix an arbitrary
irreducible component of $V$, say $V_1$. Let $\mathcal{S}$ be a
stratification of the pair $(\mathbb{CP}^{n+1},V)$. Then for a
fixed integer $1 \leq i \leq n$, the prime factors of the global
Alexander polynomial $\delta_i(t)$ of $V$ are among the prime
factors of local polynomials $\xi^s _{l}(t)$ associated to the
local Alexander modules $H_l (S^{2n-2s+1}-K^{2n-2s-1};\Gamma)$ of
 link pairs $(S^{2n-2s+1},K^{2n-2s-1})$
of components of strata $S \in \mathcal{S}$ such that: $S \subset
V_1$, $n-i \leq s=\text{dim}{S} \leq n$, and $l$ is in the range
$2n-2s-i \leq l \leq n-s$.
\end{thm}

\emph{Note.} \  The $0$-dimensional strata of $V$ may only
contribute to $\delta_n(t)$, the $1$-dimensional strata may only
contribute to $\delta_n(t)$ and $\delta_{n-1}(t)$ and so on. This
observation will play a key role in the proof of Proposition 5.1
of the next section.
\begin{proof}
We will use the Lefschetz hyperplane section theorem and induction
down on $i$. The beginning of the induction is the
characterization of the 'top' Alexander polynomial of $V$:
\emph{the prime divisors of $\delta_n(t)$ are among the prime
factors of local polynomials $\xi^s _{l}(t)$ corresponding to
strata $S \in \mathcal{S}$ with $S \subset V_1$, $0 \leq
s=\text{dim}{S} \leq n$, and $n-2s \leq l \leq n-s$.} This follows
from the following more general fact:\newline

\noindent\emph{Claim.} For any $1 \leq i \leq n$, the prime
divisors of $\delta_i(t)$ are among the prime factors of the local
polynomials $\xi^s _{l}(t)$ corresponding to strata $S \in
\mathcal{S}$ such that: $S \subset V_1$, $0 \leq s=\text{dim}{S}
\leq n$, and $i-2s \leq l \leq n-s$.

\emph{Proof of Claim.} Since $V_1$ is an irreducible component of
$V$, it acquires the induced stratification from that of $V$. By
the transversality assumption, the stratification $\mathcal{S}$ of
the pair $(\mathbb{CP}^{n+1},V)$ induces a
 stratification of the pair $(\mathbb{CP}^{n+1},V \cup H)$.

Let  $j$ and $i$ be the inclusions of $\mathbb{CP}^{n+1}-V_1$ and
respectively $V_1$ into $\mathbb{CP}^{n+1}$. For a fixed
perversity $\bar p$ we will denote the intersection complexes
$IC_{\bar p} ^{\bullet} (\mathbb{CP}^{n+1} , \mathcal{L}_H)$  by
$IC_{\bar p} ^{\bullet}$. The distinguished triangle $i_* i^! \to
id \to j_* j^* \overset{[1]}{\to}$, upon applying the
hypercohomology functor, yields the following long exact sequence:
\[
\begin{aligned}
\cdots \to \mathbb{H}^{-i}_{V_1}(\mathbb{CP}^{n+1};IC_{\bar m}
^{\bullet}) \to
IH_{i} ^{\bar m} (\mathbb{CP}^{n+1};\mathcal{L}_{H}) \to \mathbb{H}^{-i}(\mathbb{CP}^{n+1}-V_1;IC_{\bar m} ^{\bullet}) \to \\
 \to
\mathbb{H}^{-i+1}_{V_1}(\mathbb{CP}^{n+1};IC_{\bar m} ^{\bullet})
\to \cdots
\end{aligned}
\]

Note that the complex $j^*IC_{\bar m} ^{\bullet}[-n-1]$ on
$\mathbb{CP}^{n+1}-V_1$ is perverse with respect to the middle
perversity (since $j$ is the open inclusion and the functor $j^*$
is t-exact; \cite{BBD}). Therefore, by  Artin's vanishing theorem
for perverse sheaves (\cite{S}, Corollary 6.0.4) and noting that
$\mathbb{CP}^{n+1}-V_1$ is affine (\cite{Di}, (1.6.7)), we obtain:
$$\mathbb{H}^{-i} (\mathbb{CP}^{n+1}-V_1;j^*IC_{\bar m} ^{\bullet}) \cong 0 \ \ \text{for} \ \ i<n+1.$$
Therefore:
$$\mathbb{H}^{-i}_{V_1}(\mathbb{CP}^{n+1};IC_{\bar m} ^{\bullet}) \cong
IH_{i} ^{\bar m} (\mathbb{CP}^{n+1};\mathcal{L}_{H}) \ \
\text{for} \ \ i<n, $$ and $IH_{n} ^{\bar m}
(\mathbb{CP}^{n+1};\mathcal{L}_{H})$ is a quotient of
$\mathbb{H}^{-n}_{V_1}(\mathbb{CP}^{n+1};IC_{\bar m} ^{\bullet})$.

Now using the superduality isomorphism $IC_{\bar m} ^{\bullet}
\cong \mathcal{D} {IC_{\bar l} ^{\bullet}}^{op} [2n+2]$ and the
fact that the stalks over $V_1$ of the complex ${IC_{\bar l} ^{\bullet}}^{op}$ are torsion
$\Gamma$-modules, and $IC_{\bar m} ^{\bullet}|_{V_1} \cong 0$, we have the isomorphisms:
\[
\begin{aligned}
\mathbb{H}^{-i} _{V_1}(\mathbb{CP}^{n+1};IC_{\bar m} ^{\bullet})
&= \mathbb{H}^{-i} (V_1;i^!IC_{\bar m} ^{\bullet}) \\
&\cong \mathbb{H}^{-i+2n+2} (V_1;\mathcal{D}i^* {IC_{\bar l} ^{\bullet}}^{op}) \\
&\cong Hom(\mathbb{H}^{i-2n-2} (V_1;i^* {IC_{\bar l}
^{\bullet}}^{op});\Gamma) \oplus
Ext(\mathbb{H}^{i-2n-1}(V_1;i^* {IC_{\bar l} ^{\bullet}}^{op});\Gamma)\\
&\cong Ext(\mathbb{H}^{i-2n-1} (V_1;i^* {IC_{\bar l} ^{\bullet}}^{op});\Gamma)\\
&\cong \mathbb{H}^{i-2n-1} (V_1;i^* {IC_{\bar l} ^{\bullet}}^{op})\\
&\cong \mathbb{H}^{i-2n-1} (V_1;i^* {\mathcal{R}
^{\bullet}}^{op}).
\end{aligned}
\]

Therefore it suffices to study the order of the module
$\mathbb{H}^{i-2n-1} (V_1;i^* {\mathcal{R} ^{\bullet}}^{op})$, for
fixed $i \leq n$.

By the compactly supported hypercohomology long exact sequence and
induction on the strata of $V_1$, the polynomial associated to
$\mathbb{H}^{i-2n-1} (V_1;i^* {\mathcal{R} ^{\bullet}}^{op})$ will
divide the product of the polynomials associated with all the
modules $\mathbb{H}^{i-2n-1}_c (\mathcal{V};{\mathcal{R}
^{\bullet}}^{op}|_{\mathcal{V}})$, where $\mathcal{V}$ runs over
the strata of $V_1$ in the stratification of the pair
$(\mathbb{CP}^{n+1},V \cup H)$, i.e. $\mathcal{V}$ is of the form
$S \cap H$ or $S - S \cap H$, for $S \in \mathcal{S}$ and $S
\subset V_1$.

Next, we will need the following lemma:
\begin{lemma}
Let $\mathcal{V}$ be a $j$-(complex) dimensional stratum of $V_1$
(or $V$) in the stratification of the pair $(\mathbb{CP}^{n+1},V
\cup H)$. Then the prime factors of the polynomial associated to
$\mathbb{H}^{i-2n-1} _{c}(\mathcal{V};
{\mathcal{R}^{\bullet}}^{op}|_{\mathcal{V}})$ must divide one of
the polynomials $\xi^j _{l}(t) = \text{order} \{ IH_{l}^{\bar
m}(S^{2n-2j+1};\mathcal{L}) \}$, in the range $0 \leq l \leq n-j$
and $0 \leq i-l \leq 2j$, where $(S^{2n-2j+1}, K^{2n-2j-1})$ is
the link pair of $\mathcal{V}$ in ($\mathbb{CP}^{n+1},V \cup H)$.
\end{lemma}

Once the lemma is proved, the \emph{Claim} (and thus the beginning
of the induction) follows from Remark 3.3 which describes the
polynomials of link pairs of strata $S \cap H$ of $V \cap H$ in
$(\mathbb{CP}^{n+1}, V \cup H)$ in terms of the polynomials of
link pairs of strata $S \in \mathcal{S}$ of $V$ in
$(\mathbb{CP}^{n+1}, V)$, and Lemma 3.2 which relates the local
intersection Alexander polynomials of links of strata $S \in
\mathcal{S}$ to the classical local Alexander polynomials.\newline

In order to finish the proof of the theorem we use the Lefschetz
hyperplane theorem and induction down on $i$. We denote the
Alexander polynomials of $V$ by $\delta_i^V(t)$ and call
$\delta_n^V(t)$ the 'top' Alexander polynomial of $V$.

Let $1 \leq i=n-k<n$ be fixed. Consider $L \cong
\mathbb{CP}^{n-k+1}$ a generic codimension $k$ linear subspace of
$\mathbb{CP}^{n+1}$, so that $L$ is transversal to $V \cup H$.
Then $W=L \cap V$ is a $(n-k)$-dimensional, degree $d$, reduced
hypersurface in $L$, which is transversal to the hyperplane at
infinity $H \cap L$ of $L$. Moreover, by the transversality
assumption, the pair $(L,W)$ has a Whitney stratification induced
from that of the pair $(\mathbb{CP}^{n+1},V)$, with strata of the
form $\mathcal{V}=S \cap L$, for $S \in \mathcal{S}$. The local
coefficient system $\mathcal{L}_H$ defined on
$\mathcal{U}=\mathbb{CP}^{n+1}-(V \cup H)$ restricts to a
coefficient system on $\mathcal{U} \cap L$ defined by the same
representation (here we already use the Lefschetz theorem).

By applying the Lefchetz hyperplane section theorem (\cite{Di},
(1.6.5)) to $\mathcal{U}=\mathbb{CP}^{n+1}-(V \cup H)$ and its
section by $L$, we obtain the isomorphisms:
$$\pi_i(\mathcal{U} \cap L) \overset{\cong}{\to}
\pi_i(\mathcal{U}), \ \ \text{for} \ i \leq n-k.$$ and a
surjection for $i=n-k+1$. Therefore the homotopy type of
$\mathcal{U}$ is obtained from that of $\mathcal{U} \cap L$ by
adding cells of dimension $> n-k+1$. Hence the same is true for
the infinite cyclic covers $\mathcal{U}^c$ and ${(\mathcal{U} \cap
L)}^c$ of $\mathcal{U}$ and $\mathcal{U} \cap L$ respectively.
Therefore,
$$H_i({(\mathcal{U} \cap L)}^c;\mathbb{Q}) \overset{\cong}{\to} H_i(\mathcal{U}^c;\mathbb{Q}) , \ \ \text{for} \ i \leq n-k.$$
Since the maps above are induced by embeddings, these maps are
isomorphisms of $\Gamma$-modules. We conclude that
$\delta_{n-k}^W(t)=\delta_{n-k}^V(t)$.

Next, note that $\delta_{n-k}^W(t)$ is the 'top' Alexander
polynomial of $W$ as a hypersurface in $L \cong
\mathbb{CP}^{n-k+1}$, therefore by the induction hypothesis, the
prime factors of $\delta_{n-k}(t)$ are restricted to those of the
local Alexander polynomials $\xi^r_l(t)$ associated to link pairs
of strata $\mathcal{V}=S \cap L \subset W_1=V_1 \cap L$, with $0
\leq r=\text{dim}(\mathcal{V}) \leq n-k$ and $(n-k)-2r \leq l \leq
(n-k)-r$. Now, using the fact that the link pair of a stratum
$\mathcal{V}=S \cap L$ in $(L,W)$ is the same as the link pair of
$S$ in $(\mathbb{CP}^{n+1},V)$, the conclusion follows by
reindexing (replace $r$ by $s-k$, where $s=\text{dim}(S)$).
\end{proof}

\emph{Note.} The Lefschetz argument in the above proof may be
replaced by a similar argument for intersection homology modules,
using also the realization of the Alexander modules of the
hypersurface complement as intersection homology modules. More
precisely, the Lefschetz hyperplane theorem for intersection
homology  (\cite{FK2} or \cite{S}, Example 6.0.4(3)) yields the
following isomorphisms of $\Gamma$-modules:
$$IH_i^{\bar m}(L, {\mathcal{L}_H}|_L) \overset{\cong}{\to}
IH_i^{\bar m}(\mathbb{CP}^{n+1},\mathcal{L}_H), \ \ \text{for} \ i
\leq n-k.$$ On the other hand the following are isomorphisms of
$\Gamma$-modules (by Corollary 3.4): $IH_{i} ^{\bar m}
(\mathbb{CP}^{n+1};\mathcal{L}_{H}) \cong H_i (\mathcal{U} ^c;
\mathbb{Q})$ and $IH_{i} ^{\bar m} (L, {\mathcal{L}_H}|_L) \cong
H_i ((\mathcal{U} \cap L) ^c; \mathbb{Q})$.\newline

\emph{Proof of lemma 4.3.} \ For simplicity, we  let $r=i-2n-1$.
The module $\mathbb{H}_{c} ^{r}
(\mathcal{V};{\mathcal{R}^{\bullet}}^{op}|_{\mathcal{V}})$ is the
abuttement of a  spectral sequence with $E_2$ term given by:
\begin{equation*}
E_2 ^{p,q}=H^{p} _{c}(\mathcal{V}; \mathcal{H}^q
({\mathcal{R}^{\bullet}}^{op}|_{\mathcal{V}})).
\end{equation*}
Since ${\mathcal{R}^{\bullet}}^{op}$ is a constructible complex,
$\mathcal{H}^q ({\mathcal{R}^{\bullet}}^{op}|_{\mathcal{V}})$ is a
local coefficient system on $\mathcal{V}$. Therefore, by the
orientability of $\mathcal{V}$ and the Poincar\'{e} duality
isomorphism (\cite{Br}, V.9.3), $E_2 ^{p,q}$ is isomorphic to the
module $H_{2j-p}(\mathcal{V}; \mathcal{H}^q
({\mathcal{R}^{\bullet}}^{op}|_{\mathcal{V}}))$. As in Lemma 9.2
of \cite{F}, we can show that the latter is a finitely generated
module. More precisely, by deformation retracting $\mathcal{V}$ to
a closed, hence finite, subcomplex of $V_1$ (or
$\mathbb{CP}^{n+1}$), we can use simplicial homology with local
coefficients to calculate the above $E_2$ terms.

We will keep the cohomological indexing in the study of the above
spectral sequence (see, for example, \cite{F}). By the above
considerations, we may assume that $\mathcal{V}$ is a finite
simplicial complex.

$E_2 ^{p,q}$ is the $p$-th homology of a cochain complex
 $C_{c} ^*(\mathcal{V};\mathcal{H}^q ({\mathcal{R}^{\bullet}}^{op}|_{\mathcal{V}}))$
 whose $p$-th cochain group is a subgroup of
 $C^p(\mathcal{V};\mathcal{H}^q ({\mathcal{R}^{\bullet}}^{op}|_{\mathcal{V}}))$, which in turn
 is the direct sum of modules of the form
$\mathcal{H}^q ({\mathcal{R}^{\bullet}}^{op})_{x(\sigma)}$, where
$x(\sigma)$ is the barycenter of a $p$-simplex $\sigma$ of
$\mathcal{V} \subset V_1$. By the stalk calculation (\cite{B},
V.3.15) and using ${IC_{\bar m} ^{\bullet}}|_{V_1} \cong 0$,
\begin{equation*}
\mathcal{H}^q ({\mathcal{R}^{\bullet}}^{op})_{x(\sigma)} \cong
\mathcal{H} ^{q} ({IC_{\bar l} ^{\bullet}}^{op})_{x(\sigma)} \cong
\begin{cases}
0, & q > -n-j-1\\
IH_{2n+1+q}^{\bar m}(L_{x(\sigma)};\mathcal{L}),& q \leq -n-j-1
\end{cases}
\end{equation*}
(where $L_{x(\sigma)} \cong S^{2n-2j+1}$ is the link of
$\mathcal{V}$ in ($\mathbb{CP}^{n+1},V \cup H)$). Given that $E_2
^{p,q}$ is a quotient of  $C_{c} ^p (\mathcal{V};\mathcal{H}^q
({\mathcal{R}^{\bullet}}^{op}|_{\mathcal{V}}))$, we see that $E_2
^{p,q}$ is a torsion module, and a prime element ${\gamma}\in
{\Gamma}$ divides the order of $E_2 ^{p,q}$ only if it divides the
order of one of the torsion modules $IH_{2n+1+q}^{\bar
m}(L_{x(\sigma)};\mathcal{L})$. Denote by $\xi^j _{2n+1+q}(t)$ the
order of the latter module, where $j$ stands for the dimension of
the stratum.

Each $E_r ^{p,q}$ is a quotient of a submodule of $E_{r-1}
^{p,q}$, so by induction on $r$, each of them is a torsion
$\Gamma$-module whose associated polynomial has the same property
as that of $E_2$. Since the spectral sequence converges in
finitely many steps, the same property is satisfied by $E_\infty$.

By spectral sequence theory,
\begin{equation*}
E_{\infty} ^{p,q} \cong F^p \mathbb{H}^{p+q} _{c}(\mathcal{V};
{\mathcal{R}^{\bullet}}^{op}|_{\mathcal{V}})/ F^{p+1}
\mathbb{H}^{p+q} _{c}(\mathcal{V};
{\mathcal{R}^{\bullet}}^{op}|_{\mathcal{V}}),
\end{equation*}
where the modules  $F^p \mathbb{H}^{p+q} _{c}(\mathcal{V};
{\mathcal{R}^{\bullet}}^{op}|_{\mathcal{V}})$ form a descending
bounded filtration of $\mathbb{H}^{p+q} _{c}(\mathcal{V};
{\mathcal{R}^{\bullet}}^{op}|_{\mathcal{V}})$.

Now set $A^{*}=\mathbb{H}^{*} _{c}(\mathcal{V};
{\mathcal{R}^{\bullet}}^{op})$
 as a graded module which is filtered by $F^p A^{*}$ and set
$E_0^p(A^{*})=F^p A^{*}/F^{p+1}A^{*}$.  Then, for some $N$, we
have:
\begin{equation*}
{0}\subset F^N A^{*} \subset F^{N-1}A^{*} \subset \cdots \subset
F^1 A^{*}\subset F^0 A^{*}\subset F^{-1}A^{*}=A^{*}.
\end{equation*}
This yields the series of short exact sequences:

{\allowdisplaybreaks
\begin{displaymath}
\begin{diagram}[height=2em,width=2em]
0      & \rTo         &F^N A^{*} &  \rTo^{\cong}  & E_0^N(A^{*})  &\rTo  & 0                &    &   \\
0      & \rTo         &F^N A^{*} &  \rTo         & F^{N-1} A^{*} &\rTo  & E_0^{N-1}(A^{*}) &\rTo& 0 \\
       &              &          &               &  \vdots       &      &                  &    &   \\
0      & \rTo         &F^k A^{*} &  \rTo         & F^{k-1} A^{*} &\rTo  & E_0^{k-1}(A^{*}) &\rTo& 0 \\
       &              &          &               &   \vdots      &      &                  &    &   \\
0      & \rTo         &F^1 A^{*} &  \rTo         & F^{0} A^{*}   &\rTo  & E_0^{0}(A^{*})   &\rTo& 0 \\
0      & \rTo         &F^0 A^{*} &  \rTo         & A^{*}         &\rTo  & E_0^{-1}(A^{*})  &\rTo& 0 \\
\end{diagram}
\end{displaymath}
} Let us see what happens at the $r$th grade of these graded
modules. For clarity, we will indicate the grade with a
superscript following the argument. For any $p$,
\begin{align*}
E_0^p(A^{*})^r&=(F^pA^{*}/F^{p+1}A^{*})^r\\
&=F^pA^r/F^{p+1}A^r\\
&=F^pA^{p+r-p}/F^{p+1}A^{p+r-p}\\
&=E_{\infty} ^{p,r-p}.
\end{align*}
We know that each of the prime factors of the polynomial of this
module must be a prime factor of some $\xi^j _{2n+1+(r-p)}(t)$.
Further, by dimension considerations and stalk calculation, we
know that $E_\infty ^{p,r-p}$ can be non-trivial only if $0 \leq p
\leq 2j$ and $-2n-1 \leq r-p \leq -n-j-1$. Hence, as $p$ varies,
the only prime factors under consideration are those of $\xi^j
_{2n+1+(r-p)}(t)$ in this range, i.e. they are the only possible
prime factors of the $E_0^p(A^{*})^r$, collectively in $p$ (but
within the grade $r$).

By induction down the above list of short exact sequences, we
conclude that $F^N A^r=E_0^N(A^{*})^r$, and subsequently
$F^{N-1}A^r$, $F^{N-2}A^r$,$\dots$, $F^0A^r$, and $A^r$, have the
property of being torsion modules whose polynomials are products
of polynomials whose prime factors are all factors of one of the
$\xi^j _{2n+1+a}(t)$, where $a$  must be chosen in the range $0
\leq r-a \leq 2j$ and $-2n-1 \leq a \leq -n-j-1$. Since
$\mathbb{H}^{r} _{c}(\mathcal{V};
{\mathcal{R}^{\bullet}}^{op}|_{\mathcal{V}})$ is the submodule of
$A^{*}$ corresponding to the $r$th grade, it too has this
property. Using the fact that $r=i-2n-1$ and reindexing, we
conclude that the prime factors of the polynomial of
$\mathbb{H}^{i-2n-1} _{c}(\mathcal{V};
{\mathcal{R}^{\bullet}}^{op}|_{\mathcal{V}})$ must divide one of
the polynomials $\xi^j _{l}(t) = \text{order} \{ IH_{l}^{\bar
m}(S^{2n-2j+1};\mathcal{L}) \}$, in the range $0 \leq l \leq n-j$
and $0 \leq i-l \leq 2j$, where $S^{2n-2j+1}$ is the link of (a
component of) $\mathcal{V}$. $\square$

\begin{remark}\emph{\emph{Isolated singularities}}\newline
In the case of hypersurfaces with only isolated singularities,
Theorem 4.2 can be strengthen as follows.

Assume that $V$ is an $n$-dimensional reduced projective
hypersurface, transversal to the hyperplane at infinity, and
having only isolated singularities. If $n \geq 2$ this assumption
implies that $V$ is irreducible. If $n=1$, we fix an irreducible
component, say $V_1$. The only interesting global (intersection)
Alexander module is $IH_{n} ^{\bar m} (\mathbb{CP}^{n+1}
;\mathcal{L}_H) \cong H_n(\mathcal{U} ^c;\mathbb{Q})$. As in the
proof of the Theorem 4.2, the latter is a quotient of the torsion
module $\mathbb{H}^{-n-1} (V_1;{\mathcal{R}^{\bullet}}^{op})$. Let
$\Sigma_0= \text{Sing(V)} \cap V_1$ be the set of isolated
singular points of $V$ which are contained in $V_1$. Note that
$V_1$ has an induced stratification:
$$V_1 \supset (V_1 \cap H) \cup \Sigma_{0} \supset \Sigma_{0}.$$
The long exact sequence of the compactly supported hypercohomology yields:
$$\to \mathbb{H}^{-n-1} _{c} (V_1 - \Sigma_{0};{\mathcal{R}^{\bullet}}^{op}) \to
\mathbb{H}^{-n-1} (V_1;{\mathcal{R}^{\bullet}}^{op}) \to
\mathbb{H}^{-n-1} (\Sigma_{0};{\mathcal{R}^{\bullet}}^{op}) \to$$
 and by the local calculation on stalks we obtain:
\[
\begin{aligned}
\mathbb{H}^{-n-1} (\Sigma_{0};{\mathcal{R}^{\bullet}}^{op}) &\cong
\oplus_{p \in \Sigma_{0}} \mathcal{H}^{-n-1}(\mathcal{R}^{\bullet})_{p}^{op} \\ &\cong \oplus_{p \in \Sigma_{0}} IH_{n} ^{\bar
l}(S_{p} ^{2n+1};\Gamma)^{op} \\
&\overset{(1)}{\cong} \oplus_{p \in \Sigma_{0}} IH_{n} ^{\bar m}(S_{p}^{2n+1};\Gamma) \\ &\cong
\oplus_{p \in \Sigma_{0}} H_{n} (S_{p}^{2n+1}-S_{p}^{2n+1} \cap V;\Gamma),
\end{aligned}
\]
where $(S_{p}^{2n+1},S_{p}^{2n+1} \cap V)$ is the (smooth) link
pair of the singular point $p \in \Sigma_{0}$, and the module
$H_{n} (S_{p}^{2n+1}-S_{p}^{2n+1} \cap V;\Gamma)$ is the classical
(local) Alexander polynomial of the algebraic link. $(1)$ follows
from the superduality isomorphism for intersection Alexander
polynomials of link pairs (\cite{ShCa}, Corollary 3.4; \cite{F},
Theorem 5.1).\newline By Remark 3.3, Lemma 4.3 and the long exact
sequences of compactly supported hypercohomology, it can be shown
 that the modules $\mathbb{H}^{-n-1}
_{c} (V_1 - \Sigma_{0};{\mathcal{R}^{\bullet}}^{op})$ and
$\mathbb{H}^{-n} _{c} (V_1 -
\Sigma_{0};{\mathcal{R}^{\bullet}}^{op})$ are annihilated by
powers of $t-1$.

Thus we obtain the following divisibility theorem (compare
\cite{Li}, Theorem 4.3; \cite{Li4}, Theorem 4.1(1); \cite{Di2},
Corollary 6.4.16):
\begin{thm}
Let $V$ be a projective hypersurface in $\mathbb{CP}^{n+1}$ ($n
\geq 1$), which is transversal to the hyperplane at infinity, H,
and has only isolated singularities. Fix an irreducible component
of $V$, say $V_1$, and let $\Sigma_0=V_1 \cap \text{Sing}(V)$. Then $IH_{n} ^{\bar m}
(\mathbb{CP}^{n+1};\mathcal{L}_{H}) \cong H_{n} (\mathcal{U}
^c;\mathbb{Q})$ is a torsion $\Gamma$-module, whose associated
polynomial $\delta_n(t)$ divides the product $\prod_{p \in
\Sigma_0} \Delta_{p}(t) \cdot (t-1)^r$ of the local Alexander
polynomials of links of the singular points of $V$ which are
contained in $V_1$.
\end{thm}
\end{remark}

An immediate consequence of the previous theorems is the triviality
of the global polynomials $\delta_i (t)$, $1 \leq i \leq n$, if
none of the roots of the local Alexander polynomials along some
irreducible component of $V$ is a root of unity of order $d$:

\begin{example}
Suppose that $V$ is a degree $d$ reduced projective hypersurface
which is also a rational homology manifold, has no codimension $1$
singularities, and is transversal to the hyperplane at infinity.
Assume that the local monodromies of link pairs of strata
contained in some irreducible component $V_1$ of $V$ have orders
which are relatively prime to $d$ (e.g., the transversal
singularities along strata of $V_1$ are Brieskorn-type
singularities, having all exponents relatively prime to $d$).
Then, by Theorem 4.1, Theorem 4.2 and Proposition 2.1, it follows
that $\delta_i (t) \sim 1$, for $1 \leq i \leq n$ (see \cite{Mi},
theorem 9.1).
\end{example}

Further obstructions on the global Alexander modules/polynomials are provided by the relation with the
'modules/polynomials at infinity'. The following is an extension of
Theorem 4.5 of \cite{Li} or, in the case $n=1$, of Theorem 4.1(2) of \cite{Li4}.

\begin{thm}
Let $V$ be a degree $d$ reduced  hypersurface  in
$\mathbb{CP}^{n+1}$, which is transversal to the hyperplane at
infinity, $H$. Let $S_{\infty}$ be a sphere of sufficiently large
radius in $\mathbb{C}^{n+1}= \mathbb{CP}^{n+1} -H$ (or
equivalently, the boundary of a sufficiently small tubular
neighborhood of $H$ in $\mathbb{CP}^{n+1}$). Then for all $i < n$,
$$IH_{i} ^{\bar m} (\mathbb{CP}^{n+1};\mathcal{L}_{H}) \cong \mathbb{H}^{-i-1}(S_{\infty};IC^{\bullet} _{\bar m}) \cong
H_{i} (\mathcal{U}^c_{\infty};\mathbb{Q})$$ and $IH_{n} ^{\bar m}
(\mathbb{CP}^{n+1};\mathcal{L}_{H})$ is a quotient of
$\mathbb{H}^{-n-1}(S_{\infty};IC^{\bullet} _{\bar m}) \cong H_{n}
(\mathcal{U}^c_{\infty};\mathbb{Q})$, where $\mathcal{U}^c_{\infty}$ is the infinite
cyclic cover of $S_{\infty}-(V \cap S_{\infty})$ corresponding to
the linking number with $V \cap S_{\infty}$ (cf. \cite{Li}).
\end{thm}

\emph{Note.} If $V$ is an irreducible curve of degree $d$ in $\mathbb{CP}^2$, in general position at infinity, then the associated
'polynomial at infinity', i.e. the order of $H_{1} (\mathcal{U}^c_{\infty};\mathbb{Q})$, is $(t-1)(t^d-1)^{d-2}$
(see \cite{Li3}, \cite{Li4}).

\begin{proof}
Choose coordinates $(z_0: \cdots :z_{n+1})$ in the projective space such that $H=\{z_{n+1}=0 \}$ and $\mathbb{O}=(0:\cdots
:0:1)$ is the origin in $\mathbb{CP}^{n+1} -H$.  Define $$\alpha :\mathbb{CP}^{n+1} \to
\mathbb{R}_{+} \ , \ \alpha :=\frac{ \abs{z_{n+1}}^2}{\sum_{i=0}^{n+1} \abs{z_{i}}^2}$$
Note that $\alpha$ is well-defined, it is real analytic and proper, $$0 \leq \alpha \leq 1 \ , \ \alpha^{-1}(0)=H \
\text{and} \
\alpha^{-1}(1)=\mathbb{O}.$$ Since $\alpha$ has only finitely many critical values, there is $\epsilon$ sufficiently small
such that the interval $(0,\epsilon]$ contains no critical values. Set $U_{\epsilon}=\alpha^{-1}([0,\epsilon))$, a tubular
neighborhood of $H$ in $\mathbb{CP}^{n+1}$ and note that $\mathbb{CP}^{n+1}-U_{\epsilon}$ is a closed large ball of radius
$R=\frac{1-\epsilon}{\epsilon} \overset{\epsilon \to 0}{\to} \infty$ in $\mathbb{C}^{n+1}=\mathbb{CP}^{n+1} -H$.

Lemma 8.4.7(a) of \cite{KS} applied to $\alpha$ and $IC^{\bullet} _{\bar m}$, together with ${IC^{\bullet} _{\bar m}}|_{H}
\cong 0$, yield:
$$\mathbb{H}^{\ast}(U_{\epsilon};IC^{\bullet} _{\bar m}) \cong \mathbb{H}^{\ast}(H;IC^{\bullet} _{\bar m}) \cong 0$$
and therefore, by the hypercohomology long exact sequence, we obtain the isomorphism:
$$\mathbb{H}^{\ast}(\mathbb{CP}^{n+1};IC^{\bullet} _{\bar m}) \cong
\mathbb{H}^{\ast} _{\mathbb{CP}^{n+1}-U_{\epsilon}} (\mathbb{CP}^{n+1};IC^{\bullet} _{\bar m})$$
Note that, for $i: \mathbb{CP}^{n+1}-U_{\epsilon} \hookrightarrow \mathbb{CP}^{n+1}$ the inclusion,
\begin{align*}
\mathbb{H}^{\ast} _{\mathbb{CP}^{n+1}-U_{\epsilon}} (\mathbb{CP}^{n+1};IC^{\bullet} _{\bar m}) & =
\mathbb{H}^{\ast} (\mathbb{CP}^{n+1}-U_{\epsilon}; i^!IC^{\bullet} _{\bar m}) \cong
\mathbb{H}^{\ast} (\mathbb{CP}^{n+1};i_*i^!IC^{\bullet} _{\bar m})\\ &\cong
\mathbb{H}^{\ast} (\mathbb{CP}^{n+1};i_!i^!IC^{\bullet} _{\bar m}) \overset{def}{=}
\mathbb{H}^{\ast} (\mathbb{CP}^{n+1},U_{\epsilon};IC^{\bullet} _{\bar m})\\  &\cong
\mathbb{H}^{\ast} (\mathbb{C}^{n+1},U_{\epsilon}-H;IC^{\bullet} _{\bar m}),
\end{align*}
where the last isomorphism is the excision of $H$ (see for example \cite{Ma}, \S 1; \cite{Di}, Remark 2.4.2(ii)).

If $k$ is the open inclusion of the affine space in $\mathbb{CP}^{n+1}$, then $k^*IC_{\bar m} ^{\bullet}[-n-1]$
is perverse with respect to the middle perversity (since $k$ is the open
inclusion and the functor $k^*$ is t-exact). Therefore, by  Artin's
vanishing theorem for perverse sheaves (\cite{S}, Corollary 6.0.4), we obtain:
$$\mathbb{H}^{-i} (\mathbb{C}^{n+1};k^*\mathcal{IC}^{\bullet} _{\bar m}) \cong 0 \ \ \text{for} \ \ i<n+1.$$
The above vanishing and the long exact sequence of the pair $(\mathbb{C}^{n+1},U_{\epsilon}-H)$ yield the
isomorphisms:
$$\mathbb{H}^{-i} _{\mathbb{CP}^{n+1}-U_{\epsilon}} (\mathbb{CP}^{n+1};IC^{\bullet} _{\bar m}) \cong
\mathbb{H}^{-i-1} (U_{\epsilon}-H;IC^{\bullet} _{\bar m}) \ \ \ \text{if} \ \ i<n$$
and
$$\mathbb{H}^{-n-1} (U_{\epsilon}-H;IC^{\bullet} _{\bar m}) \to
\mathbb{H}^{-n} _{\mathbb{CP}^{n+1}-U_{\epsilon}} (\mathbb{CP}^{n+1};IC^{\bullet} _{\bar m})
\ \ \text{is an epimorphism.}$$
Note that $U_{\epsilon}-H=\alpha^{-1}((0,\epsilon))$ and by Lemma 8.4.7(c)
of \cite{KS} we obtain the isomorphism:
$$\mathbb{H}^{\ast} (U_{\epsilon}-H;IC^{\bullet} _{\bar m}) \cong \mathbb{H}^{\ast} (S_{\infty};IC^{\bullet} _{\bar m})$$
where $S_{\infty}=\alpha^{-1}(\epsilon ')$, $0 < \epsilon ' < \epsilon $.

Next, using the fact that
${IC^{\bullet} _{\bar m}}|_{{V \cup H}} \cong 0$, we obtain a sequence of isomorphisms as follows: for $i \leq n$,
\begin{align*}
\mathbb{H}^{-i-1}(S_{\infty};IC^{\bullet} _{\bar m})
&=\mathbb{H}^{-i-1} _{c}(S_{\infty}-(V \cap S_{\infty});IC^{\bullet} _{\bar m})\\
&=\mathbb{H}^{-i-1} _{c}(S_{\infty}-(V \cap S_{\infty});\mathcal{L}|_{{S_{\infty}-V}}[2n+2])\\
&=H^{-i+2n+1} _{c}(S_{\infty}-(V \cap S_{\infty});\mathcal{L})\\
&\overset{(1)}{\cong} H_{i} (S_{\infty}-(V \cap S_{\infty});\mathcal{L})\\
&\cong H_{i} (\mathcal{U}^c_{\infty};\mathbb{Q}),
\end{align*}
where $\mathcal{L}$ is given on $S_{\infty}-(V \cap S_{\infty})$
by the linking number with $V \cap S_{\infty}$, $(1)$ is the
Poincar\'{e} duality isomorphism (\cite{Br}, Theorem V.9.3), and
$\mathcal{U}^c_{\infty}$ is the infinite cyclic cover of $S_{\infty}-(V \cap
S_{\infty})$ corresponding to the linking number with $V \cap
S_{\infty}$ (cf. \cite{Li}).
\end{proof}

\begin{remark} Subsequently, A. Libgober has found a simpler proof of Theorem 4.8,
using a purely topological argument based on the Lefschetz
theorem. As a corollary to Theorem 4.8 it follows readily (cf.
\cite{Li6}) that the  Alexander modules of the hypersurface
complement are semi-simple, thus generalizing Libgober's result
for the case of hypersurfaces with isolated singularities (see
\cite{Li}, Corollary 4.8). The details will be given below.
\end{remark}

\begin{prop}
Let $V \subset \mathbb{CP}^{n+1}$ be a degree $d$ reduced hypersurface which is transversal to the hyperplane at infinity, $H$.
Then for each $i \leq n$, the Alexander module $H_{i} (\mathcal{U}^c;\mathbb{C})$ is a semi-simple  $\mathbb{C}[t,t^{-1}]$-module which is
annihilated by $t^d-1$.
\end{prop}

\begin{proof}
By Theorem 4.8, it suffices to prove this fact for the  modules
'at infinity' $H_{i} (\mathcal{U}^c_{\infty};\mathbb{C})$, $i \leq
n$.

Note that since $V$ is transversal to $H$, the space
$S_{\infty}-(V \cap S_{\infty})$ is a circle fibration over $H - V \cap H$ which is homotopy equivalent to the
complement in $\mathbb{C}^{n+1}$
to the affine cone over the  projective hypersurface $V\cap H$. Let $\{h=0\}$ be the polynomial defining $V\cap H$ in $H$.
Then the infinite cyclic cover
$\mathcal{U}^c_{\infty}$ of $S_{\infty}-(V \cap S_{\infty})$ is homotopy equivalent to the Milnor fiber $\{h=1\}$ of the
(homogeneous) hypersurface singularity at the origin defined by $h$
and, in particular, $H_{i} (\mathcal{U}^c_{\infty};\mathbb{C})$ ($i \leq n$) is a torsion finitely generated $\mathbb{C}[t,t^{-1}]$-module.
Since the monodromy on the Milnor fiber $\{h=1\}$ has finite order $d$ (given by multiplication by roots of unity),
it also follows that the  modules at infinity are semi-simple torsion modules,  annihilated by $t^d-1$ (see \cite{Ku}).
\end{proof}

\emph{Note.} The above proposition supplies alternative proofs to Corollary 4.8 and Theorem 4.1.


\section{On the Milnor fiber of a projective arrangement of hypersurfaces}
In this section, we apply the preceding results to the case of a
hypersurface $V \subset \mathbb{CP}^{n+1}$, which is a projective
cone over a reduced hypersurface $Y \subset \mathbb{CP}^n$. As an
application to Theorem 4.4, we obtain restrictions on the
eigenvalues of the monodromy operators associated to the Milnor
fiber of the hypersurface arrangement defined by $Y$ in
$\mathbb{CP}^n$.

Let $f : \mathbb{C}^{n+1} \to \mathbb{C}$ be a homogeneous polynomial of degree $d>1$,
and let $Y=\{f=0\}$ be the projective hypersurface in $\mathbb{CP}^n$ defined by $f$.
Assume that the polynomial $f$ is square-free and let $f=f_1\cdots f_s$ be the decomposition
of $f$ as a product of irreducible factors. Then $Y_i=\{f_i=0\}$ are precisely the irreducible
components of the hypersurface $Y$, and we refer to this situation by saying that we have a
\emph{hypersurface arrangement} $\mathcal{A}=(Y_i)_{i=1,s}$ in $\mathbb{CP}^n$.

The \emph{Milnor fiber of the arrangement} $\mathcal{A}$ is defined as the fiber
$F=f^{-1}(1)$ of the global Milnor fibration $f:\mathcal{U} \to \mathbb{C}^*$ of the
(homogeneous) polynomial $f$; here $\mathcal{U}:=\mathbb{C}^{n+1}-f^{-1}(0)$ is the
complement of the central arrangement $A=f^{-1}(0)$ in $\mathbb{C}^{n+1}$, the cone
on $\mathcal{A}$. $F$ has as characteristic homeomorphism $h:F \to F$ the mapping
given by $h(x)=\tau \cdot x$ with $\tau=exp(2\pi i/d)$. This formula shows that
$h^d=id$ and hence the induced morphisms $h_q:H_q(F) \to H_q(F)$ at the homology
level are all diagonalizable over $\mathbb{C}$, with eigenvalues among the $d$-th roots of unity. Denote by $P_q(t)$
the characteristic polynomial of the monodromy operator $h_q$.

Note that the Milnor fiber $F$ is homotopy equivalent to the infinite cyclic cover
$\mathcal{U}^c$
of $\mathcal{U}$, corresponding to the homomorphism
$\mathbb{Z}^s=H_1(\mathcal{U}) \to \mathbb{Z}$
sending a meridian generator about a component of $A$ to the positive generator of
$\mathbb{Z}$. With this identification, the monodromy homeomorphism $h$ corresponds precisely to
a generator of the group of covering transformations (see \cite{Di}, p. 106-107).

It's easy to see that $V \subset \mathbb{CP}^{n+1}$, the
projective cone on $Y$, is in general position at infinity, where
we identify the hyperplane at infinity, $H$, with the projective
space on which $Y$ is defined as a hypersurface. Denote the
irreducible components of $V$ by $V_i$, $i=1,\cdots,s$, each of
which is the projective cone over the corresponding component of
$Y$. Theorem 4.2 when applied to $F \simeq \mathcal{U}^c$ and to
the hypersurface $V$, provides obstructions on the eigenvalues of
the monodromy operators associated to the Milnor fiber $F$. More
precisely, we obtain the following result concerning the prime
divisors of the polynomials $P_q(t)$, for $q \leq n-1$ (compare
\cite{Li100}, Theorem 3.1):
\begin{prop}
Let $Y=(Y_i)_{i=1,s}$ be a hypersurface arrangement in
$\mathbb{CP}^n$, and fix an arbitrary component, say $Y_1$.  Let
$F$ be the  Milnor fibre of the arrangement. Fix a Whitney
stratification of the pair $(\mathbb{CP}^n,Y)$ and denote by
$\mathcal{Y}$ the set of (open) singular strata. Then for $q \leq
n-1$, a prime $\gamma \in \Gamma$ divides the characteristic
polynomial $P_q(t)$ of the monodromy operator $h_q$ only if
$\gamma$ divides one of the polynomials $\xi^s _{l}(t)$ associated
to the local Alexander modules
$H_{l}(S^{2n-2s-1}-K^{2n-2s-3};\Gamma)$ corresponding to
 link pairs $(S^{2n-2s-1},K^{2n-2s-3})$
of components of strata $\mathcal{V} \in \mathcal{Y}$ of complex
dimension $s$ with $\mathcal{V} \subset Y_1$, such that: $n-q-1
\leq s \leq n-1$ and $2(n-1)-2s-q \leq l \leq n-s-1$.
\end{prop}
\begin{proof}
There is an identification $P_q(t) \sim \delta_q(t)$, where
$\delta_q(t)$ is the global Alexander polynomial of the
hypersurface $V$, i.e. the order of the torsion module
 $H_q(\mathcal{U}^c;\mathbb{Q}) \cong IH_q ^{\bar m}(\mathbb{CP}^{n+1};\mathcal{L}_H)$.
 We consider a topological stratification $\mathcal{S}$ on $V$
 induced by that of $Y$, having the cone point as a zero-dimensional stratum.
 From Theorem 4.2 we
 recall that, for $q \leq n-1$, the local polynomials of the zero-dimensional strata of $V_1$ do
 not contribute
 to the prime
factors of the global polynomial $\delta_q(t)$. Notice that link
pairs of strata $S$ of $V_1$ in $(\mathbb{CP}^{n+1},V)$ (with
$\text{dim}(S) \geq 1$) are the same as the link pairs of strata
of $Y_1=V_1 \cap H$ in $(H=\mathbb{CP}^n, V \cap H =Y)$. The
desired conclusion follows from Theorem 4.2 by reindexing.
\end{proof}
\emph{Note.} The polynomials $P_i(t)$, $i=0,\cdots,n$ are related by the formula
(see \cite{Di}, (4.1.21) or \cite{Di2}, (6.1.10)):
$$\prod_{q=0}^{n}P_q(t)^{(-1)^{q+1}}=(1-t^d)^{-\chi(F)/d}$$
where $\chi(F)$ is the Euler characteristic of the Milnor fiber. Therefore,
it suffices to compute only the polynomials $P_0(t),\cdots,P_{n-1}(t)$ and the Euler characteristic
of $F$.\newline

If $Y \subset \mathbb{CP}^n$ has only isolated singularities, the
proof of the previous proposition can be strengthened to obtain
the following result, similar to \cite{Di}, (6.3.29) or
\cite{Di2}, Corollary 6.4.16:

\begin{prop}
With the above notations, if $Y$ has only isolated singularities,
then the polynomial $P_{n-1}$ divides (up to a power of $t-1$) the
product of the local Alexander polynomials associated to the
singular points of $Y$ contained in $Y_1$.
\end{prop}

A direct consequence of Proposition 5.1 is the following:
\begin{cor}
If $\lambda \neq 1$ is a $d$-th root of unity such that $\lambda$
is not an eigenvalue of any of the local monodromies corresponding
to link pairs of singular strata of $Y_1$ in a stratification of
the pair $(\mathbb{CP}^n,Y)$, then $\lambda$ is not an eigenvalue
of the monodromy operators acting on $H_q(F)$ for $q \leq n-1$.
\end{cor}

Using the fact that normal crossing divisor germs have trivial monodromy operators (\cite{Di},
(5.2.21.ii); \cite{Di2}, (6.1.8.i)), we also obtain the following (compare \cite{COD}, Corollary 16):
\begin{cor} Let $\mathcal{A}=(Y_i)_{i=1,s}$ be a hypersurface arrangement  in
$\mathbb{CP}^n$ and fix one irreducible component, say $Y_1$.
Assume that $\bigcup_{i=1,s} Y_i$ is a normal crossing divisor at
any point $x \in Y_1$. Then the monodromy action on
$H_q(F;\mathbb{Q})$ is trivial for $q \leq n-1$.
\end{cor}


\section{Examples}
We will now show, by explicit calculations on examples, how to
combine Theorems $4.1$ and $4.2$ in order to obtain information on
the Alexander modules of a hypersurface.

We start with few remarks on the local Alexander polynomials of
link pairs of strata of $(\mathbb{CP}^{n+1},V)$. Let $s$ be the
complex dimension of a (component of a) stratum $\mathcal{V} \in
\mathcal{S}$, and consider $(S^{2n-2s+1},K)
\overset{not}{=}(S,K)$, the link pair in $(\mathbb{CP}^{n+1},V)$
of a point $x \in \mathcal{V}$. This is in general a singular
algebraic link, obtained by  intersecting $(\mathbb{CP}^{n+1},V)$
with a small sphere centered at $x$, in a submanifold of
$\mathbb{CP}^{n+1}$ of dimension $n-s+1$, which meets the
$s$-dimensional stratum transversally at $x$. We define, as usual,
a local system on the link complement, $S-K$, with stalk $\Gamma$
and action of the fundamental group given by $\alpha \mapsto
t^{lk(\alpha,K)}$. The classical Alexander modules of the link
pair are defined as $H_* (S-K,\Gamma) \cong H_*
(\widetilde{S-K},\mathbb{Q})$, where $\widetilde{S-K}$ is the
infinite cyclic covering of the link complement defined
geometrically by the total linking number with $K$, i.e. the
covering associated to the kernel of the epimorphism $\pi_1 (S -
K) \overset{\text{lk}}{\to} \mathbb{Z}$, which maps the meridian
generator loops around components of $K$ to $1$. The
$\Gamma$-module structure on $H_* (\widetilde{S-K},\mathbb{Q})$ is
induced by the action of the covering transformations. Note that
$H_* (S-K,\Gamma)$ is torsion $\Gamma$-module since algebraic
knots are of finite type (\cite{ShCa}, Proposition 2.2).

For every algebraic link $(S,K)$ as above, there is an associated
Milnor fibration (\cite{Mi}): ${^{s}F} \hookrightarrow S-K \to
S^1$, where $^{s}F$ is the Milnor fiber. Whenever we refer to
objects associated to an s-dimensional stratum, we will use the
superscript '$^{s}$'. It is known that $^{s}F$ has the homotopy
type of a $(n-s)$-dimensional CW complex (\cite{Mi}, Theorem 5.1).
We regard $H_* (^{s}F;\mathbb{Q})$ as a $\Gamma$-module, with the
multiplication defined by the rule: $t \times a = {^{s}h}_*(a)$,
where $^{s}h : {^{s}F} \to {^{s}F}$ is the monodromy homeomorphism
of the Milnor fibration. Note that the inclusion $^{s}F
\hookrightarrow S-K$ is, up to homotopy equivalence, the infinite
cyclic covering $\widetilde{S-K}$  of the link complement, defined
by the total linking number with $K$ (more precisely,
$\widetilde{S-K}$ is homeomorphic to $^{s}F \times \mathbb{R}$;
note that $^{s}F$ is connected since we work with reduced
singularities).
 With this
identification, the monodromy homeomorphism of $^{s}F$ corresponds
precisely to a generator of the group of covering transformations.
It follows that the classical Alexander polynomials of the link
pair can be identified with the characteristic polynomials
${^{s}\delta} _* (t)=\text{det}(tI-{^{s}h}_*)$ of the monodromy
operators ${^{s}h}_*$.

Note that, in general, the Milnor fibre associated to an algebraic link has a certain degree of
connectivity. In the case of an isolated singularity, the Milnor fibre is homotopy equivalent to a join of spheres of its middle dimension
(\cite{Mi}). Results on
the homotopy type of the Milnor fibre of a non-isolated hypersurface singularity and homology calculations can be found,
for example,  in \cite{Oka} and \cite{Ra}.\newline

\noindent\textbf{Example 6.1: }\emph{One-dimensional singular locus}\newline
Let $V$ be the trifold in $\mathbb{CP}^4 = \{(x:y:z:t:v)\}$, defined by the polynomial:
$y^{2}z+x^3+tx^2+v^3=0$. The singular locus of $V$ is the projective line
$\text{Sing}(V)=\{(0:0:z:t:0); z,t \in \mathbb{C}\}$. We
let $\{t=0\}$ be the hyperplane $H$ at infinity. Then $V \cap H$ is the surface in $\mathbb{CP}^3 = \{(x:y:z:v)\}$
defined by the equation
$y^{2}z+x^3+v^3=0$, having the point $(0:0:1:0)$ as its singular set. Thus,
$\text{Sing}(V \cap H)=\text{Sing}(V) \cap H$. Let $X$ be the affine part of $V$, i.e., defined by the polynomial
$y^{2}z+x^3+x^2+v^3=0$. Then $\text{Sing}(X)=\{(0,0,z,0)\} \cong \{(0:0:z:1:0)\}=\text{Sing}(V) \cap X$ is the $z$-axis of
$\mathbb{C}^4 =\{(x,y,z,v)\}$,
and it's clear that the origin $(0,0,0,0)=(0:0:0:1:0)$ looks different than any other point on the $z$-axis:
the tangent cone at the point $(0,0,\lambda,0)$ is represented by two planes for $\lambda \neq 0$
and degenerates to a double plane for $\lambda =0$. Therefore we give the pair $(\mathbb{CP}^4,V)$
the following Whitney stratification:
$$\mathbb{CP}^4 \supset V \supset \text{Sing}(V) \supset (0:0:0:1:0)$$
It's clear that $V$ is transversal to the hyperplane at infinity.

In our example ($n=3$, $k=1$) we are interested in describing the prime factors of the global Alexander polynomials
$\delta_2(t)$ and $\delta_3(t)$ (note that $\delta_1(t) \sim 1$, as $n-k \geq 2$; cf. \cite{Li}).
In order to describe the local Alexander polynomials of link pairs of singular strata of $V$,
 we will use the results of \cite{Oka} and \cite{Ra}.

The link pair of the top stratum of $V$ is $(S^1, \emptyset)$, and the only prime factor that may contribute to
the global Alexander polynomials is $t-1$, the order of $H_0 (S^1,\Gamma)$.

Next, the link of the stratum $\text{Sing}(V)-\{(0:0:0:1:0)\}$ is the algebraic knot in a $5$-sphere $S^5 \subset
 \mathbb{C}^3$ given by the intersection of the affine variety $\{y^2+x^3+v^3=0\}$ in $\mathbb{C}^3 =\{(x,y,v)\}$
 (where $t=0$, $z=1$) with a small sphere about
 the origin $(0,0,0)$. (To see this, choose the  hyperplane $V(t)=\{(x:y:z:0:v)\}$ which is transversal
 to the singular set $V(x,y,v)$, and consider an affine neighborhood of their intersection $(0:0:1:0)$ in $V(t) \cong
 \mathbb{CP}^3$.) The polynomial $y^2+x^3+v^3$ is weighted homogeneous of Brieskorn type, hence (\cite{Mi}, \cite{Oka})
 the associated Milnor fibre is simply-connected, homotopy equivalent to $\{2 \ \text{points} \} * \{3 \ \text{points} \}
 * \{3 \ \text{points} \}$,
 and the characteristic
 polynomial of the monodromy operator acting on $H_2 (F; \mathbb{Q})$ is $(t+1)^{2} (t^2-t+1)$ .

 Finally, the link pair of the zero-dimensional stratum, $\{(0:0:0:1:0)\}$, (the origin in the affine space $\{t=1\}$),
 is the
 algebraic knot in a $7$-sphere, obtained by intersecting the affine variety $y^{2}z+x^3+x^2+v^3=0$ in $\mathbb{C}^4 =
 \{(x,y,z,v)\}$ with a
 small sphere about the origin. Since we work in a neighborhood of the origin, by an analytic change of coordinates,
 this is
 the same as the link pair of the origin in the variety $y^{2}z+x^2+v^3=0$.
 Therefore the Milnor fiber of the associated Milnor
 fibration is the join of $\{x^2=1\}$, $\{y^2z=1\}$ and $\{v^3=1\}$, i.e., it is homotopy equivalent to
 $S^2 * \{3 \ \text{points} \}$, i.e., $S^3 \vee S^3$
 (\cite{Mi}, \cite{Oka}). Moreover, denoting by $\delta_f$ the characteristic polynomials of monodromy of the weighted homogeneous
 polynomial $f$, we obtain (\cite{Oka}, Theorem $6$): $\delta_{x^2+v^3+y^2z}(t) = \delta_{x^2+v^3+yz}(t) = \delta_{x^2+v^3}(t) =
 t^2-t+1$.

 Note that the above links, $K^3 \subset S^5$ and $K^5 \subset S^7$, are rational homology spheres since none of the
 characteristic polynomials of their associated Milnor fibers has the trivial eigenvalue $1$ (see \cite{Ra}, Proposition
 3.6). Therefore $V$, and hence $V \cap H$, are rational homology manifolds (see the discussion preceding Proposition 2.1).
 Then, by Proposition 2.1, $t-1$ cannot be a prime factor of the global Alexander polynomials of $V$. Also note that the
 local Alexander polynomials of links of the singular strata of $V$ have prime divisors none of which divides $t^3-1$, thus,
 by Theorem 4.1, they cannot appear among the prime divisors of $\delta_2(t)$ and $\delta_3(t)$.

 Altogether, we conclude that $\delta_0(t) \sim t-1$, $\delta_1(t) \sim 1$, $\delta_2(t) \sim 1$ and $\delta_3(t) \sim
 1$.\newline

\emph{Note.} The above example can be easily generalized to provide hypersurfaces of any dimension,
with a one-dimensional singular locus and trivial global Alexander polynomials. This can be done by adding cubes of new
variables to the polynomial $y^{2}z+x^3+tx^2+v^3$.\newline

\noindent\textbf{Example 6.2: }\emph{An isolated singularity}\newline
Let $V$ be the surface in $\mathbb{CP}^3$, defined by the polynomial:
$x^3+y^3+z^3+tz^2=0$. The singular locus of $V$ is a point $\text{Sing}(V)=\{(0:0:0:1)\}$. We
let $\{t=0\}$ be the hyperplane $H$ at infinity and note that $V \cap H$ is a nonsingular hypersurface in $H$ (defined by
the zeros of polynomial $x^3+y^3+z^3$). Hence $H$ is transversal on $V$ (in the stratified sense). The only nontrivial
global Alexander polynomial (besides $\delta_0(t) \sim t-1$) may be
$\delta_2(t)$, and its prime divisors are either $t-1$ or prime divisors of the local Alexander polynomial of the link
pair of the isolated singular point of $V$. The link pair of $\{(0:0:0:1)\}$ in $(\mathbb{CP}^3,V)$ is the
 algebraic knot in a $5$-sphere, obtained by intersecting the affine variety $x^3+y^3+z^3+z^2=0$ in $\mathbb{C}^3$ with a
 small sphere about the origin. Since we work in a neighborhood of the origin, by an analytic change of coordinates, this is
 the same as the link pair of the origin in the variety $x^3+y^3+z^2=0$.
The polynomial $x^3+y^3+z^2$ is weighted homogeneous of Brieskorn type, hence (\cite{Mi}, \cite{Oka}) the characteristic
 polynomial of the monodromy homeomorphism of the associated Milnor fibration is $(t+1)^{2} (t^2-t+1)$. Note that none of
 the roots of the latter polynomial is a root of unity of order $3$, hence none of its prime factors may appear as a prime
 factor of $\delta_2(t)$ (by Theorem 4.1). Moreover, the link of the isolated singular point is a rational homology sphere,
 as the monodromy operator of the associated Milnor fibration has no trivial eigenvalue (\cite{Di}, Theorem 3.4.10 (A)).
 Therefore, by Proposition 2.1, $t-1$ is not a prime factor of $\delta_2(t)$.
 Altogether, we conclude that $\delta_2(t) \sim 1$ and the Alexander module is trivial.

The same answer is obtained by using Corollary 4.9 of \cite{Li}. Indeed, we have an isomorphism of
 $\Gamma$-modules: $H_{2}  (\mathbb{CP}^{3} - V \cup H;\Gamma) \cong H_2 (M_f ; \mathbb{Q})$, where $M_f$ is the
 Milnor fiber at the origin, as a non-isolated hypersurface singularity in $\mathbb{C}^4$,
 of the polynomial $f : \mathbb{C}^4 \to \mathbb{C}$, $f(x,y,z,t)=
 x^3+y^3+z^3+tz^2$. The module structure on $H_2 (M_f ; \mathbb{Q})$ is given by the action on the monodromy operator.
 By a linear change of coordinates, we can work instead with the polynomial $x^3+y^3+tz^2$. Note that the Milnor fiber of
 the latter is $\{ 3 \ \text{points} \} * \{ 3 \ \text{points} \} * S^1$. Hence the formula of the homology of a join (\cite{Di}, (3.3.20))
 shows that the module $H_2 (M_f ; \mathbb{Q})$ is trivial.\newline

In the above examples, our theorems are used to show the triviality of the global
Alexander modules. But these global objects are not always trivial.\newline

\noindent\textbf{Example 6.3: }\emph{Manifold singularity}\newline
Consider the hypersurface $V$ in $\mathbb{CP}^{n+1}$ defined by the zeros of the polynomial: $f(z_0,
\cdots, z_n)=(z_0)^2+(z_1)^2+\cdots +(z_{n-k})^2$. Assume that $n-k \geq 2$ is even.
The singular set $\Sigma=V(z_0, \cdots ,z_{n-k}) \cong
\mathbb{CP}^{k}$ is non-singular. Choose a generic hyperplane, $H$, for example $H=\{z_n = 0 \}$.
The link of $\Sigma$ is the algebraic knot in a $(2n-2k+1)$-sphere $S_{\epsilon} ^{2n-2k+1}$ given by
the intersection of the affine variety $(z_0)^2+(z_1)^2+\cdots +(z_{n-k})^2=0$ in $\mathbb{C}^{n-k+1}$
with a small sphere about the origin. As $n-k$ is even, the link of the singularity (in the sense of \cite{Mi})
is a rational homology sphere and
the associate Alexander polynomial of the knot complement is $t-(-1)^{n-k+1}=t+1$. Hence the prime factors of the
intersection Alexander polynomials of the hypersurface are either $t+1$ or $t-1$. However, since the links of singular
strata are rational homology spheres, we conclude (by using Proposition 2.1 and \cite{Li}, Corollary 4.9) that
$\delta_{n-k}(t) \sim t+1$ and all the
$\delta_{j}(t) \ , \ n-k <j \leq n$, are multiples of $t+1$. Also note that in this case, $\delta_{j}(t) \sim 1$ for
$1 \leq j \leq n-k-1$.\newline

 Sometimes it is possible to calculate explicitly the global Alexander polynomials, as we will
 see in the next example:\newline

\noindent\textbf{Example 6.4: }
Let $V$ be the surface in $\mathbb{CP}^3$, defined by the homogeneous polynomial of degree $d$:
$f(x,y,z,t)=x^{d-1}z+xt^{d-1}+y^{d}+xyt^{d-2}$, $d \geq 3$. The singular locus of $V$ is a point:
$\text{Sing}(V)=\{(0:0:1:0)\}$. We fix a generic hyperplane, $H$. Then the (intersection) Alexander modules of $V$
are defined and the only 'non-trivial' Alexander polynomial of the hypersurface is $\delta_2(t)$. Note that, by
 Corollary 4.9 of
 \cite{Li},  there is an isomorphism of
 $\Gamma$-modules: $H_{2}(\mathbb{CP}^{3}- V \cup H;\Gamma) \cong H_2 (M_f ; \mathbb{Q})$, where $M_f$ is the
 Milnor fiber at the origin, as a non-isolated hypersurface singularity in $\mathbb{C}^4$,
 of the polynomial $f : \mathbb{C}^4 \to \mathbb{C}$, $f(x,y,z,t)= x^{d-1}z+xt^{d-1}+y^{d}+xyt^{d-2}$, and where
 the module structure on $H_2 (M_f ; \mathbb{Q})$ is given by the action on the monodromy operator.
Moreover, by \cite{Di}, Example 4.1.26,  the characteristic polynomial of the latter is $t^{d-1}+t^{d-2}+ \cdots
+1$. Therefore, $\delta_2(t)= t^{d-1}+t^{d-2}+ \cdots +1$.


\providecommand{\bysame}{\leavevmode\hbox to3em{\hrulefill}\thinspace}

\bigskip \noindent
\small\textsc{Department of Mathematics, \\
University of Pennsylvania, \\
209 S 33rd St, \\
Philadelphia, PA 19104-6395, USA} \\
{\em e-mail: }{\tt lmaxim@math.upenn.edu}

\end{document}